\documentclass[a4book, 10pt]{amsart}
\title[Flat M\"obius strips]{
   Flat M\"obius strips of given isotopy type in $\R^3$ \\
whose centerlines are geodesics or \\  lines of curvature}
\date{2007 August 12.}
\usepackage[dvips]{graphicx} 
\usepackage{verbatim,enumerate}
\usepackage{amssymb}
\usepackage{times}
\usepackage{amsthm}
\theoremstyle{plain}
 \newtheorem{theorem}{Theorem}[section]
 \newtheorem*{theorem*}{Theorem}

 \newtheorem*{lemma*}{Lemma}
 \newtheorem{proposition}[theorem]{Proposition}
 
 \newtheorem*{fact*}{Fact}

 \newtheorem{lemma}[theorem]{Lemma}
 \newtheorem{corollary}[theorem]{Corollary}
\theoremstyle{remark}
 \newtheorem{definition}[theorem]{Definition}
 \newtheorem{remark}[theorem]{Remark}
 \newtheorem*{remark*}{Remark}
 \newtheorem*{acknowledgements}{Acknowledgements}
 \newtheorem{example}[theorem]{Example}
\numberwithin{equation}{section}

\newcommand{\pt}[1]{\mathsf{#1}}

\newcommand{\Z}{\boldsymbol{Z}}
\newcommand{\R}{\boldsymbol{R}}

\renewcommand{\qed}{q.e.d.}

\renewcommand{\phi}{\varphi}
\renewcommand{\epsilon}{\varepsilon}

\newcommand{\op}{\operatorname}

\newcommand{\mb}[1]{{\mathbf #1}}

\newcommand{\pmt}[1]{{\begin{pmatrix} #1  \end{pmatrix}}}

\newcommand{\dy}{\displaystyle}

\author{Yasuhiro Kurono }
\address[Kurono]{%
   Department of Mathematics, Graduate School of Science,
   Osaka University,
   Toyonaka, Osaka 560-0043,
   Japan
}

\author{Masaaki Umehara}
\address[Umehara]{%
   Department of Mathematics, Graduate School of Science,
   Osaka University,
   Toyonaka, Osaka 560-0043,
   Japan 
}
\email{umehara@math.sci.osaka-u.ac.jp} 

\thanks{
The second author were supported by Grant-in-Aid for 
Scientific Research (A) No.~19204005 
from the Japan Society for the Promotion of Science.
}

\begin{document}
\begin{abstract}
 We construct flat M\"obius strips of arbitrary 
isotopy types,
whose centerlines are geodesics or 
lines of curvature. 
\end{abstract}
\maketitle

\section*{Introduction}
Let $\phi(p,t)=(-p,-t)$ ($|t|<1,\,\,p\in S^1$) be an involution
on $S^1\times (-1,1)$, where
$$
S^1:=\{p\in \R^2\,;\, |p|=1\}.
$$
We denote by $\Bbb M$ the quotient of
$S^1\times (-1,1)$ by the map $\phi$, and
denote by 
$$
\pi:S^1\times (-1,1)\to \Bbb M:=S^1\times (-1,1)/\sim_\phi
$$
the canonical projection.
An embedding $f:\Bbb M\to \R^3$ is called a {\it M\"obius strip}
and the restriction of $f$ on the line $\{(p,t)\in S^1\times (-1,1)\,;\,t=0\}$
is called the {\it centerline} of the M\"obius strip.
A M\"obius strip $f$ is called {\it rectifying} (or {\it geodesic}) if 
the centerline is a geodesic.
On the other hand, a M\"obius strip is called a {\it M\"obius developable} if 
it is a ruled surface and its Gaussian curvature
vanishes identically.
It should be remarked that constructing a concrete example
of M\"obius developable is not so easy, and  classical
such examples are given in Wunderlich \cite{W},
Chiconne-Kalton \cite{CK}, Schwarz \cite{S1,S2},and
Randrup-R\o gen \cite{RR}.
 
A M\"obius developable $f$ is called {\it principal} 
(or {\it orthogonal}) if the centerline 
is orthogonal to the asymptotic line.
On the complement of the set of umbilics on $\Bbb M$,
the centerline of the principal developable $f$
consists of a line of curvature. 
It should be remarked that {\it any M\"obius 
developable has at least one umbilcs} 
(See Corollary 3.5 in \cite{MU} and also Proposition \ref{prop:umbilics}
in Section 1.) In this paper, we shall prove the following two theorems:

\medskip
\noindent
{\bf Theorem A.} {\it
There exists a principal real-analytic M\"obius developable
which is isotopic to a given M\"obius strip.}

It should be remarked that the first example of unknotted principal 
real-analytic M\"obius developable was
given in \cite{CK}.

\medskip
\noindent
{\bf Theorem B.} {\it
There exists a rectifying real-analytic M\"obius developable
which is isotopic to a given M\"obius strip.}

\medskip
When we ignore the property of
the centerline, the existence of a $C^\infty$ M\"obius developable
with a given isotopy type has been shown: 
In fact, Chicone and Kalton  showed (in 1984 see \cite{CK}) 
that the existence of M\"obius developable whose center line is
an arbitrary given generic space curves.
After that, R\o gen \cite{R} showed that any embedded surfaces with boundary
in $\R^3$ can be isotopic to flat surfaces.

If we expand a flat M\"obius developable 
into their asymptotic directions, then
we get a flat surface whose asymptotic lines
are all complete, and such a surface may have
singular points in general.
In \cite{MU}, the global properties of such 
surfaces are investigated.

As a point of view from paper-handicraft,
we know experimentally  
the existence of a developable M\"obius strip 
which can be given as
an isometric deformation of
a rectangular domain on a plane. 
Such a M\"obius strip
must be rectifying, since the property that
the centerline is a geodesic is preserved by
the isometric deformation.
On the other hand, any rectifying 
M\"obius developable can be obtained by
an isometric deformation of
a rectangular domain on a plane (See Proposition 
\ref{prop:rectangular}). 
Thus, Theorem B implies that
{\it one can construct a developable M\"obius strip 
of given isotopy type via a rectangular ribbon.}

 \section{Preliminaries}\label{sec:prelim}

Let $I:=[a,b]$ be a closed interval, and
$\gamma(t)$ ($a\le t\le b$)  a regular space curve.
Then the function
$$
\kappa(t):=\frac{|\dot \gamma(t)\times \ddot \gamma(t)|}{|\dot \gamma(t)|^3}$$
is called the {\it curvature function} of $\gamma$.
A point where $\kappa(t)$ vanishes
is called an {\it inflection point} of $\gamma$, where
$\dot \gamma=d\gamma/dt$.

Let $\xi(t)$  be a vector field in $\R^3$ along 
the curve $\gamma(t)$.
We set
$$
F_{\gamma,\xi}(t,u):=\gamma(t)+u \xi(t) 
\qquad (t\in I, \,\, |u|<\varepsilon),
$$
where $\varepsilon$ is a sufficiently small positive constant.
Then $F_{\gamma,\xi}$  is called a {\it ruled strip} if it satisfies
$$
\dot\gamma(t)\times \xi(t)\ne 0,
$$
where $\times$ is the vector product in $\R^3$.
In this case, $F_{\gamma,\xi}$ gives an immersion for
sufficiently small $\epsilon$.
Moreover, if it satisfies
\begin{equation}\label{eq:flat}
\op{det}(\dot\gamma(t),\xi(t),\dot \xi(t))=0\qquad
(a\le t\le b),
\end{equation}
then $F_{\gamma,\xi}$ is called a {\it developable strip}.
In fact, it is well-known that \eqref{eq:flat} is equivalent to
the condition that the Gaussian curvature of 
$F_{\gamma,\xi}$ vanishes identically.

\begin{definition}\label{def:principal}
Let $F_{\gamma,\xi}$ be a developable strip.
Then it is called {\it principal} or {\it orthogonal} if it satisfies 
\begin{equation}\label{eq:orthogonal}
\xi(t)\cdot \dot \gamma(t)=0\qquad
(a\le t\le b),
\end{equation}
where $\cdot$ means the canonical inner product in $\R^3$.
In fact, the condition \eqref{eq:orthogonal} 
is the orthogonality of the
centerline with respect to the asymptotic direction.
If $\gamma(t)$ is not an umbilic, the centerline is
a line of curvature near $\gamma(t)$.
\end{definition}

The following assertion can be proved directly:

\begin{proposition}\label{prop:principal}
Let $\gamma$ be a regular space curve, and $\xi(t)$
a vector field along $\gamma(t)$ such that
\begin{equation}\label{eq:principal}
\xi(t)\cdot \dot \gamma(t)=0,\qquad
\dot \xi(t)\times \dot \gamma(t)=0\qquad (a\le t\le b).
\end{equation}
Then $F_{\gamma,\xi}$ gives a principal developable strip.
\end{proposition}

\begin{remark}
One can prove that any 
principal developable strip is given in this manner.
\end{remark}

\begin{remark}\label{rmk:principal}
The condition \eqref{eq:principal}
means that $\xi(t)$ is parallel with respect to the normal 
connection. In particular, the length $|\xi(t)|$ is constant
along $\gamma$. 
When $\gamma$ does not admit inflection points,
the torsion function of $\gamma$ is defined by
$$
\tau(t):=
\frac{\op{det}(\dot \gamma(t),\ddot \gamma(t),\dddot \gamma(t))}
{|\dot \gamma(t)\times \ddot \gamma(t)|^2}.
$$
We now take $t$ to be the arclength parameter.
Then, as pointed out in \cite{CK},
$$
P_0(t):=\left(\sin \int_a^t \tau(s)ds\right)
 \mb n(t)+\left (\cos \int_a^t \tau(s)ds \right) \mb b(t)
$$
gives a parallel vector field on the normal bundle $T^\perp_\gamma$
of $\gamma$, that is, $\dot P_0(t)$ is proportional to $\dot \gamma(t)$.
(Here $\mb n(t)$ and $\mb b(t)$ are the principal normal
vector field and the bi-normal vector field of $\gamma(t)$,
respectively.)
It can be easily checked that
any parallel vector field satisfying
\eqref{eq:principal} is expressed by 
$$
P(t):=(\cos \delta)P_0(t)+(\sin \delta)
\biggl(\dot\gamma(t)\times
P_0(t)\biggr),
$$
for a suitable constant $\delta\in [a,b)$.
Let $\xi(t)$ ($a\le t\le b$) 
be a non-vanishing normal vector field along $\gamma$,
that is, it satisfies $\xi(t)\cdot \dot \gamma=0$.
Let $\alpha(t)$ be the leftward angle of
$\xi(t)$ from $P(t)$.
We set
$$
\op{Tw}_{\gamma}(\xi):=\alpha(b)-\alpha(a)
$$
which is called the {\it total twist} of $\xi$
along $\gamma$, and is equal to
the total change of angles of $\xi(t)$
towards the clockwise direction with respect to
$P_0(t)$.
When $|\xi(t)|=1$,
it is well known that the following identity holds:
\begin{equation}\label{eq:twist}
\op{Tw}_{\gamma}(\xi)=\frac{1}{2\pi}\int_a^b \op{det}
(\dot \gamma(t), \xi(t),\dot \xi(t))\,dt.
\end{equation}
\end{remark}

\begin{definition}\label{def:rectifying}
Let $F_{\gamma,\xi}$ be a developable strip.
Then it is called {\it rectifying} (or {\it geodesic}) 
if it satisfies 
$$
\dot\xi(t)\cdot \ddot \gamma(t)=0\qquad
(a\le t \le b),
$$
where $\cdot$ means the canonical inner product in $\R^3$.
\end{definition}

First, we give a trivial (but important) example:

\begin{example}\label{ex:cylinder}(The cylindrical strips)
Let $\gamma(t)={}^t(x(t), y(t), 0)$ be a regular curve 
which lies entirely in the $xy$-plane. Then the cylinder
$$
F(t,u):=\gamma(t)+\pmt{0 \\ 0\\ u}
$$
over $\gamma$ gives
a developable strip {\it which is principal and rectifying 
at the same time.}
It is called a {\it cylindrical strip}.
\end{example}

Again, we return to the general setting:
Let $\gamma(t)$ ($a\le t \le b$) be a regular space curve.
If the torsion function $\tau(t)$ of $\gamma(t)$ does not vanish,
then the rectifying developable over $\gamma$ is
uniquely determined as follows:
We set
$$
D(t)=\frac{\tau(t)}{\kappa(t)}\mb t(t) + \mb b(t),
$$
which is called the
{\it normalized Darboux vector field} (cf. Izumiya-Takeuchi
\cite{IT}), 
where $\mb t(t):=\dot\gamma(t)/|\dot\gamma(t)|$.
The original Darboux vector 
field is equal to $\mb n(t)\times \dot{\mb n}(t)$,
which is proportional to $D(t)$,
where $\mb t(t),\mb n(t),\mb b(t)$
are the unit tangent vector, the unit principal normal
vector and the unit bi-normal vector, respectively.

Then one can easily get the following assertion:

\begin{proposition}\label{prop:rectifying}
Let $\gamma(t)$ be a regular space curve without
inflection points, and $D(t)$
the normalized Darboux vector field along $\gamma$.
Then $F_{\gamma,D}$ gives a rectifying developable strip.
\end{proposition}

\begin{remark}
One can prove that any
rectifying developable strip is given in this manner.
\end{remark}

Let $F_{\gamma,\xi}$ be a developable strip over 
a regular space curve $\gamma(t)$ ($ a\le t\le b$).
If it holds that
$$
\gamma^{(n)}(a)=\gamma^{(n)}(b)\qquad
(n=0,1,2,...)
$$
then $\gamma$ gives a smooth closed curve, where
$\gamma^{(n)}(t):=d^n\gamma/dt^n$.
Moreover, if 
\begin{equation}\label{eq:xi}
\xi^{(n)}(a)=-\xi^{(n)}(b)\qquad
(n=0,1,2,...)
\end{equation}
holds, then $F_{\gamma,\xi}$ gives a M\"obius developable
as defined in Introduction.
We denote by the boundary of $F_{\gamma,\xi}$ by $B_\gamma$.
The half of the linking number
$$
\op{Mtn}(F_{\gamma,\xi}):=\frac12 \op{Link}(\gamma,B_\gamma)
$$
is called the {\it M\"obius twisting number}, which takes values in
$\pm \frac12,\pm \frac32,\pm \frac52,\cdots$
(cf. \cite[Definition 3]{R}).
Here $\op{Mtn}(F_{\gamma,\xi})=(2n+1)/2$ implies that the 
strip is $(2n+1)\pi$-twisted into clockwise direction. 
Let $\mb c$ be a unit vector in $\R^3$ and suppose that the
projection of the centerline $\gamma$ into the plane
$P_\mb c$ perpendicular to $\mb c$ gives a generic plane curve.
Then we get a knot diagram of $\gamma$ on the plane $P_\mb c$,
and its writhe
$
\op{Wr}_{\mb c}(\gamma)
$
is defined, which is the total sum of the sign of crossings on
the knot diagram.
Then the following identity is well-known:
\begin{equation}\label{eq:index}
\op{Mtn}(F_{\gamma,\xi})=
-\op{Tw}_{\gamma}(\xi^{\perp})
+\op{Tw}_{\gamma}(\mb c^{\perp})+\op{Wr}_{\mb c}(\gamma),
\end{equation} 
where $\xi^{\perp}$ and $\mb c^\perp$ mean
the projection of vectors $\xi(t),\mb c$ into
the normal plane $T^\perp_\gamma$ at $\gamma(t)$.

Here, we shall recall the following result:
\begin{proposition}(\cite[Corollary 3.5]{MU})\label{prop:umbilics}
Any M\"obius developable admits at least one umbilical point.
\end{proposition}
\begin{proof}
For the sake of convenience, we shall give here a proof.
Let $\gamma(t)$ ($a\le t\le b$) be the centerline of
the M\"obius developable.
We may regard $\gamma(t)$ is a $c$-periodic
regular space curve ($c=b-a$), that is
$$
\gamma(t+c)=\gamma(t)\qquad (t\in \R).
$$ 
Then the M\"obius developable can be written as
$$
F(t,u)=\gamma(t)+u \xi(t) \qquad (|u| <\epsilon),
$$ 
where $\xi(t)$ is a unit vector field along $\gamma$ such that
\begin{equation}\label{odd:xi}
\xi(t+c)=-\xi(t)\qquad (t\in \R).
\end{equation} 
Let $\nu(t)$ be the unit normal vector field 
of $F(t,u)$, which depends 
only on $t$.
Suppose that $f$ has no umbilics.
Then we can take a local curvature line coordinate
$(x,y)$. Then by the Weingarten formula, we have
\begin{equation}\label{eq:Weingarten}
\nu_x=-\lambda_1 f_x,\qquad \nu_y=-\lambda_2 f_y,
\end{equation}
where $\lambda_1,\lambda_2$ are principal curvatures.
Without loss of generality, we may assume that $\lambda_1=0$.
Then $f_x(t,u)$ is proportional to $\xi(t)$.
Since $\lambda_1=0$,  \eqref{eq:Weingarten} yields that 
\begin{equation}\label{eq:new}
\dot \nu(t)=\nu_x \dot x+\nu_y \dot y=\nu_y \dot y=\dot y\lambda_2 f_y,
\end{equation} 
namely, $\dot \nu$ is proportional to the non-zero principal direction
$f_y$.
Since the two principal directions are orthogonal,
$\xi(t)$ must be orthogonal to $\nu(t)$ 
and $\dot \nu(t)$.
Since we have assumed that $f$ has no umbilical point,
$\nu(t)\times \dot\nu(t)$ never vanishes for all $t$.
Thus, we can write
\begin{equation}\label{eq:cross}
\xi(t)=a(t) \nu(t)\times \dot \nu(t),
\end{equation}
where $a(t)$ is a smooth function.
Since $f$ is non-orientable,
$\nu(t)$ is odd-periodic (that is $\nu(t+c)=-\nu(t)$).
In particular, $
\nu(t)\times \dot \nu(t)
$ must be $c$-periodic, that is
\begin{equation}\label{cross:nu}
\nu(t+c)\times \dot \nu(t+c)=\nu(t)\times \dot \nu(t)
\qquad (t\in \R).
\end{equation}
By \eqref{odd:xi}, \eqref{eq:cross}
and \eqref{cross:nu}, the function $a(t)$ must satisfy
the property $a(t+c)=-a(t)$.
In particular, there exists $t_0\in [a,b)$ such that
$a(t_0)=0$.
Thus we have $\xi(t_0)=0$, which contradicts that $\xi$ is
a unit vector field.
\end{proof}

Now, we would like to recall a method for constructing 
real analytic rectifying M\"obius developables from \cite{RR}.
We now assume that $\gamma(t)$ $(a\le t\le b)$
gives an embedded closed real analytic
regular space curve, which 
has no inflection points on $(a,b)$.
Since {\it a rectifying M\"obius developable must have
at least one inflection point} (See \cite{RR}),
$t=a$ must be the inflection point of $\gamma$.
Let $D(t)$ ($a< t<b$) be the normalized Darboux vector
field of $\gamma$.
Then $F_{\gamma,D}$ gives a rectifying 
M\"obius developable if and only if
$\xi:=D$ satisfies \eqref{eq:xi}, which 
reduces to  the following  Lemma \ref{lem:RR}:
The first non-vanishing non-zero coefficient vector 
$\mb c(\ne 0)$ of
the expansion of
$
\dot\gamma(t)\times \ddot\gamma(t)
$
at $t=a$ satisfies
$$
\dot\gamma(t)\times \ddot\gamma(t)=
\mb c (t-a)^N+\mbox{higher order terms},
$$
where the integer $N(\ge 1)$ is called
the {\it order of the inflection point}
and the point $t=a$ is called a {\it generic
inflection point}.
(The number $N$ is independent of the choice of the
parameter $t$ of the curve.)
Next we set
$$
\Delta(t):=\op{det}(\dot\gamma(t),\ddot\gamma(t),\dddot\gamma(t)),
$$
which is the numerator in the definition of the torsion function.
(See Remark \ref{rmk:principal}.)
Then there exists a nonzero constant $c_1$ such that
$$
\Delta(t)=
c_1 (t-a)^M+\mbox{higher order terms},
$$
where the integer $M(\ge 1)$ is called
the {\it order of torsion} at $t=0$.
The following assertion is very useful:

\begin{lemma}\label{lem:RR}{\rm (Randrup-R\o gen [RR])}
Let $\gamma(t)$ $(a\le t \le b)$
be a closed regular space such that $t=a$ is 
an inflection point, and there are no other inflection point
on $(a,b)$.
Then the  normalized Darboux vector field $D(t)$ can be 
smoothly extended as a $C^\infty$-vector field
around $t=a$ if and only if $M/N\ge 3$.
In this case, $F_{\gamma,D}$ defines a rectifying developable.
Moreover, if $N$ is odd, $F_{\gamma,D}$ is non-orientable.
\end{lemma}

As a corollary, we prove the following
assertion, which will play an important role in Section 3.

\begin{corollary}\label{cor:generic}
Suppose that the inflection point at $t=a$
is generic $($that is, $N=1)$.
Then $F_{\gamma,D}$ gives a rectifying
M\"obius developable if and only if
$$
\op{det}(\dot\gamma(t),\gamma^{(3)}(t),
\gamma^{(4)}(t)) 
$$
vanishes at $t=a$.
\end{corollary}

\begin{proof}
Since $t=a$ is an inflection point, we have 
$\ddot\gamma(a)=0$.
In particular,
$$
\dot \Delta(t)=
\op{det}(\dot\gamma(t),\ddot\gamma(t),\gamma^{(4)}(t)) 
$$
vanishes at $t=a$. On the other hand, we have
$$
\ddot \Delta(a)=
\op{det}(\dot\gamma(a),\gamma^{(3)}(a),
\gamma^{(4)}(a)), 
$$
which vanishes if and only if $M\ge 3$.
\end{proof}

Here, we give a few examples. 

\begin{figure}[ht]
\begin{center}
 \includegraphics[width=4cm]{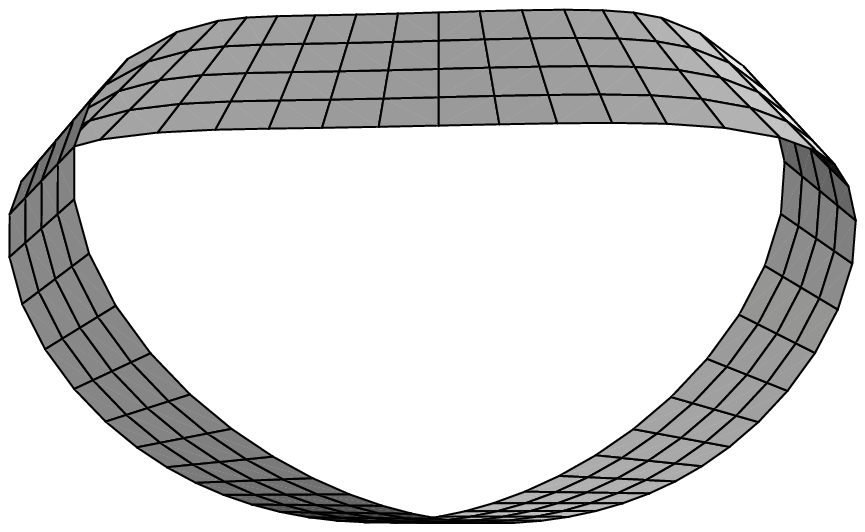}
\qquad \quad
 \includegraphics[width=4cm]{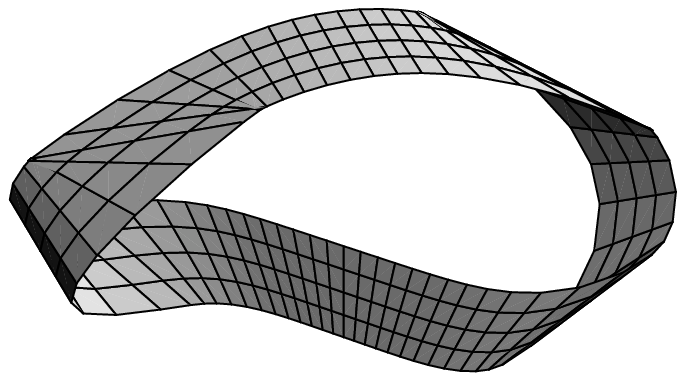} 
\end{center}
\caption{The M\"obius strips given in Examples 
\ref{ex:Wunderlich} and \ref{ex:kurono}.}
\label{Ex:bands}
\end{figure}

\begin{example}\label{ex:Wunderlich}(Wunderlich \cite{W})
Consider a regular space curve
$$
\gamma(t)=\frac1{\delta(t)}
\pmt{ 3t+2t^3+t^5\\4t+2t^3\\ -24/5 }
\qquad (t\in \R),
$$
where  $\delta(t)=9+4t^2+4t^4+t^6$. Then 
$\gamma(t)$ has no inflection point for $t\in \R$.
Moreover, it can be smoothly extended as an
embedding in $\R^3$.
In fact, 
$$
\gamma(1/s)=\frac1{\hat\delta(s)}
\pmt{3 s^5+2 s^3+s\\ 4 s^5+2s^3\\ -24 s^6/5}
\qquad (\hat \delta(s):=9s^6+4s^4+4s^2+1)
$$
is smooth at $s=0$.
This point $s=0$ is a generic inflection point 
with $N=1$ and $M=4$, and the induced 
rectifying
M\"obius developable is unknotted and of
M\"obius twisting number $1/2$.
See Figure \ref{Ex:bands} left.
\end{example}

Next, we shall give a new example of
a rectifying M\"obius developable 
whose centerline has a non-generic inflection point.

\begin{example}\label{ex:kurono}
Consider a regular space curve
$$
\gamma(t)=\frac{1}{\delta(t)}
\pmt{t^9+t^7+t^5+t^3+t\\t^5+t^3+t\\1}
\qquad (t\in \R),
$$
where $\delta(t):=t^{10}+t^8+t^6+t^4+t^2+1$.
Like as in the previous example, $\gamma(1/s)$ is
also real analytic at $s=0$ and $\gamma$ gives 
an embedded closed space curve in $\R^3$.
Moreover, $s=0$ is an inflection point
with $(N,M)=(3,10)$, that is, it is not
a generic inflection point.
By Lemma \ref{lem:RR}, the curve induces a real analytic
M\"obius developable which is unknotted and of
M\"obius twisting number $1/2$.
See Figure \ref{Ex:bands} right.
\end{example}

Randrup-R\o gen \cite{RR} gave other  examples
of rectifying M\"obius developable via Fourier polynomials.

\medskip
As pointed out in the introduction,
any M\"obius developable constructed from an isometric deformation
of rectangular domain on a plane is rectifying. 
Conversely, we can prove the following, namely,
any M\"obius developable is an isometric deformation
of rectangular domain on a plane:
\begin{proposition}\label{prop:rectangular}
Let $F=F_{\gamma,D}:[a,b]\times (\epsilon,\epsilon)$ 
be an $($embedded$)$ rectifying
M\"obius developable.
Then there exists a point $t_0\in [a,b)$ such that
the asymptotic direction $\xi(t_0)$ at $f(t_0,0)$
is perpendicular to $\dot \gamma(t_0)$.
In particular, the image $\{f(t,u)\in \R^3\,;\, t\ne t_0\}$
contains a subset which is
isometric to a rectangular domain in a plane.
\end{proposition}
\begin{proof}
Since $f$ is non-orientable,
the unit asymptotic vector filed $\xi(t)$ is
odd-periodic, that is, $\xi(a)=-\xi(b)$.
Then we have
$$
\xi(0)\cdot \dot \gamma(0)=-\xi(\pi)\cdot \dot \gamma(\pi),
$$
which implies that the function 
$t\mapsto \xi(t)\cdot \dot \gamma(t)$
changes sign on $[a,b)$.
By the intermediate vale theorem, there exists
a point $t_0\in [0,\pi)$ such that
$$
\xi(t_0)\cdot \dot \gamma(t_0)=0,
$$
which proves the assertion.
\end{proof}

\section{A $C^\infty$ M\"obius developable 
of a given isotopy-type}

In this section, we construct a rectifying $C^\infty$
M\"obius developable 
of a given isotopy-type.
To accomplish this, we prepare a special
kind of developable strip as follows:

\begin{figure}[ht]
\begin{center}
 \includegraphics[width=2.7cm]{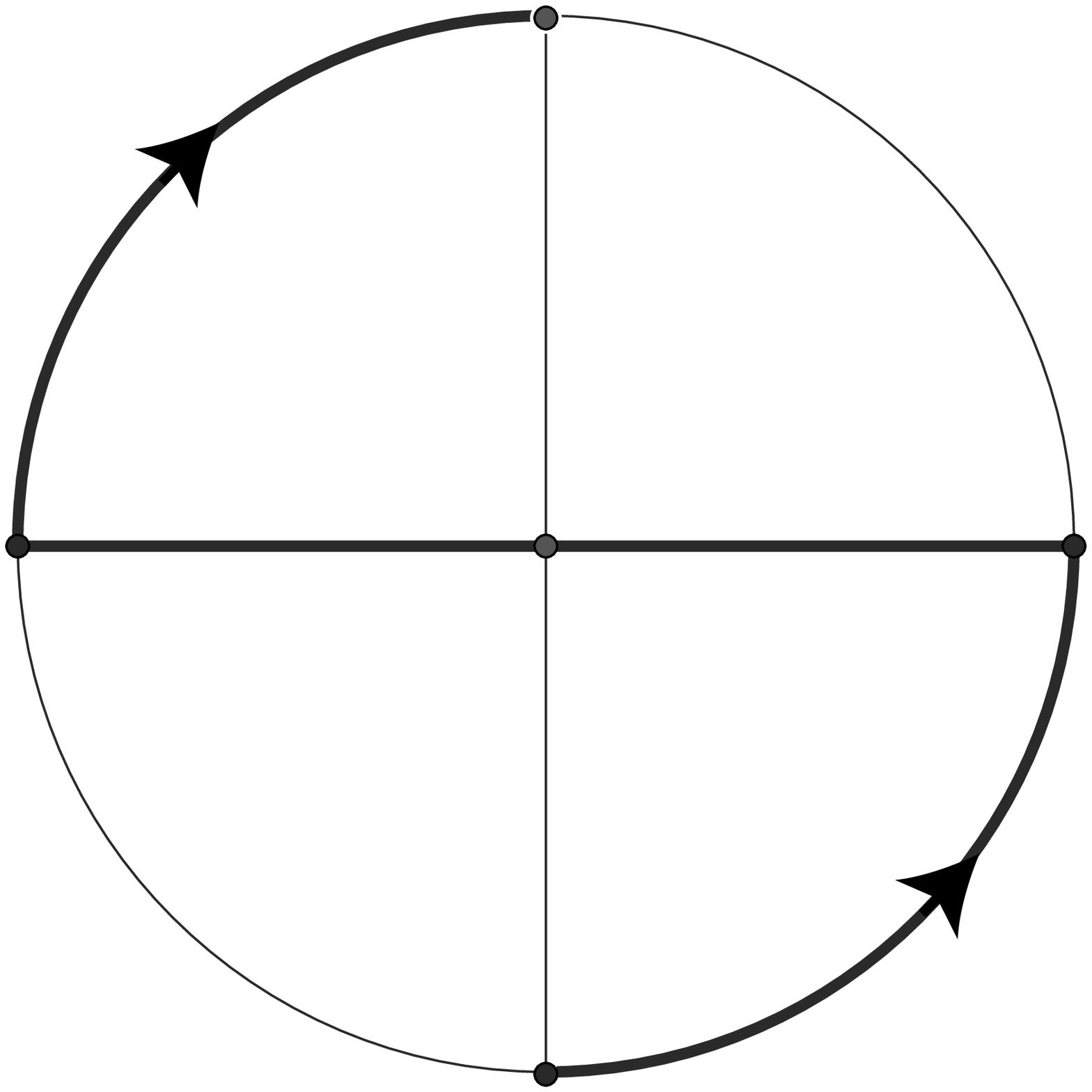}
\,\,\,\raisebox{1.5cm}{$\Longrightarrow$}
\,\,\,
 \includegraphics[width=2.7cm]{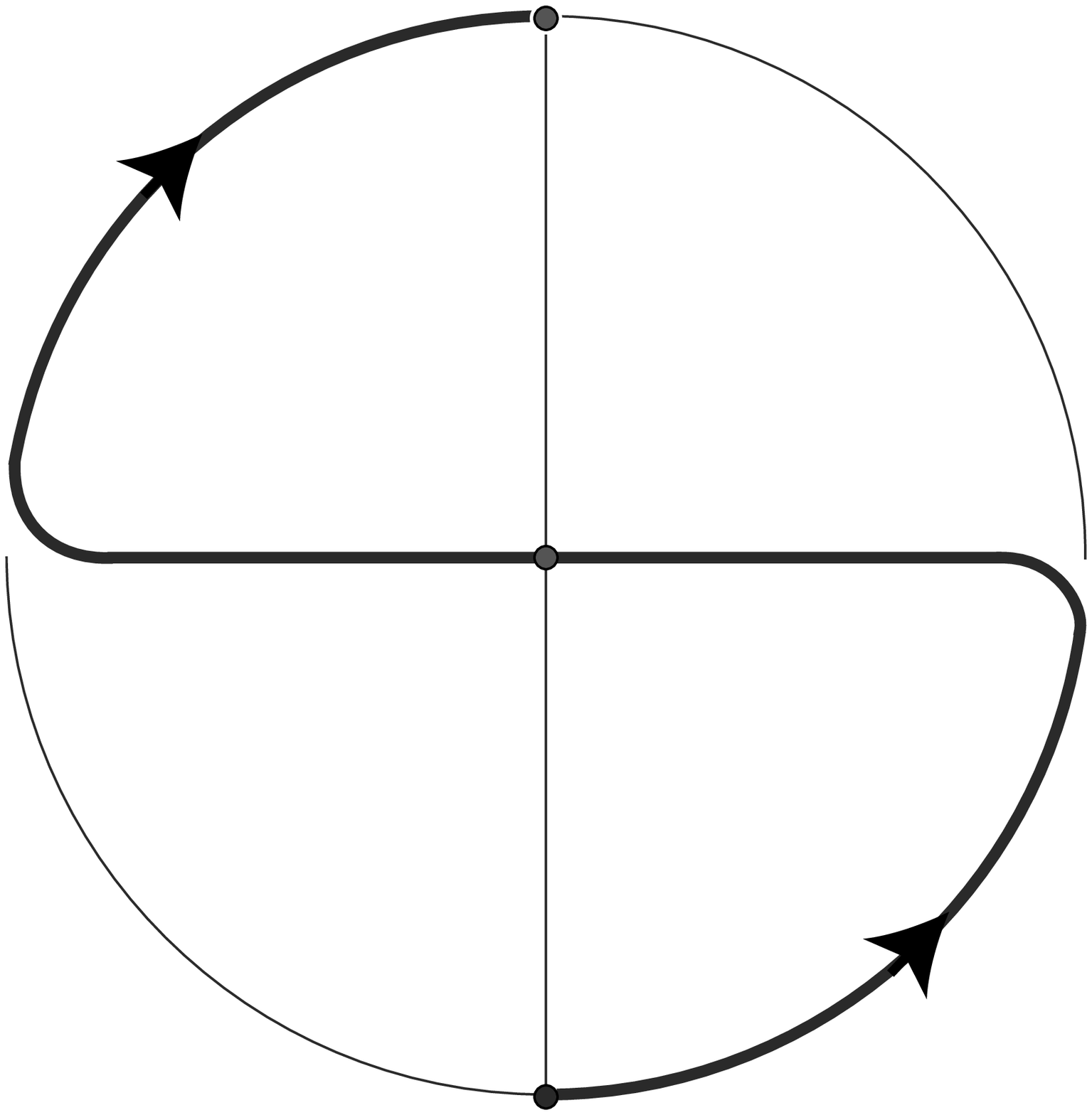} 
\caption{The original arc $\sigma$ (left) 
and $\hat \sigma$ (right)}
\end{center}
\label{twist}
\end{figure}

\medskip
\noindent
{\bf (The twisting arcs)}
Let $S^2_+$ (resp. $S^2_-$) be an upper (resp. a lower) open 
hemisphere of the unit
sphere, and let
$$
\pi_{\pm}:S^2_\pm \to \bar \Delta^2:
=\biggl\{(x,y)\in \R^2\,;\, x^2+y^2\le 1 \biggr\}
$$
be two canonical orthogonal projections.
Consider an oriented (piece-wise smooth) planar
curve $\sigma$ on the closed unit disc
$\bar \Delta^2$ as in Figure \ref{twist}.
Let $\hat \sigma$ be a $C^\infty$-regular curve rounding corner as in
right-hand side of Figure \ref{twist}.
Then the oriented space curves as the inverse images
$$
\tilde \sigma_+:=\pi_+^{-1}(\hat \sigma),\qquad 
\tilde \sigma_-:=\pi_-^{-1}(\hat \sigma)
$$
are called the {\it leftward twisting arc} or
the {\it rightward twisting arc}, respectively.
(See Figure \ref{twist}.)

\begin{figure}
[htb] 
\begin{center}
 \includegraphics[height=2.8cm]{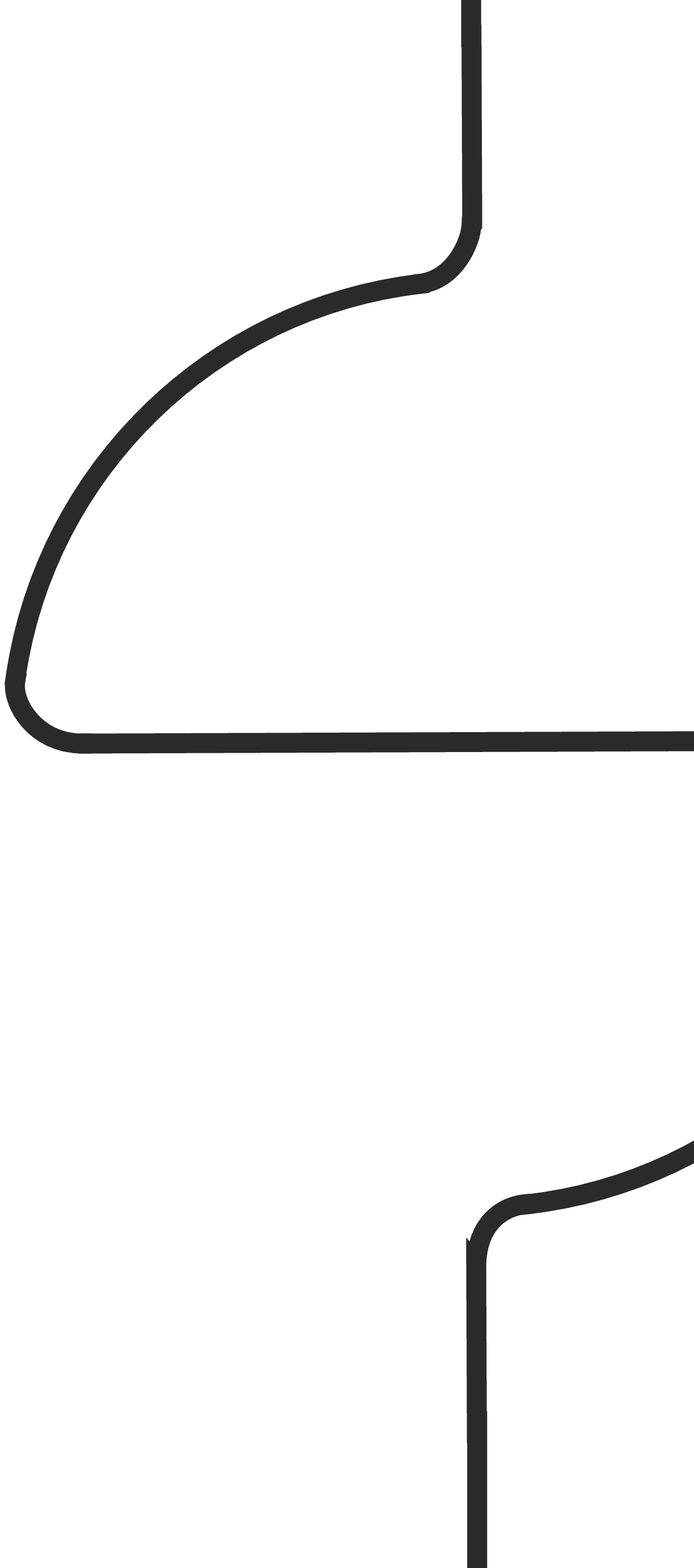}
\,\,\,\raisebox{1.5cm}{$\Longrightarrow$}
\,\,\,
 \includegraphics[height=2.8cm]{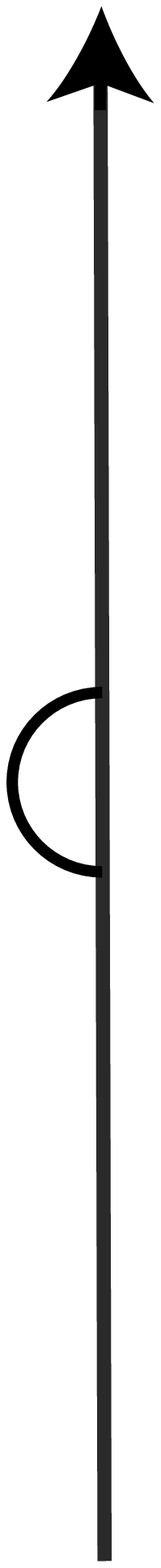} 
\hskip 2.2cm
 \includegraphics[height=2.8cm]{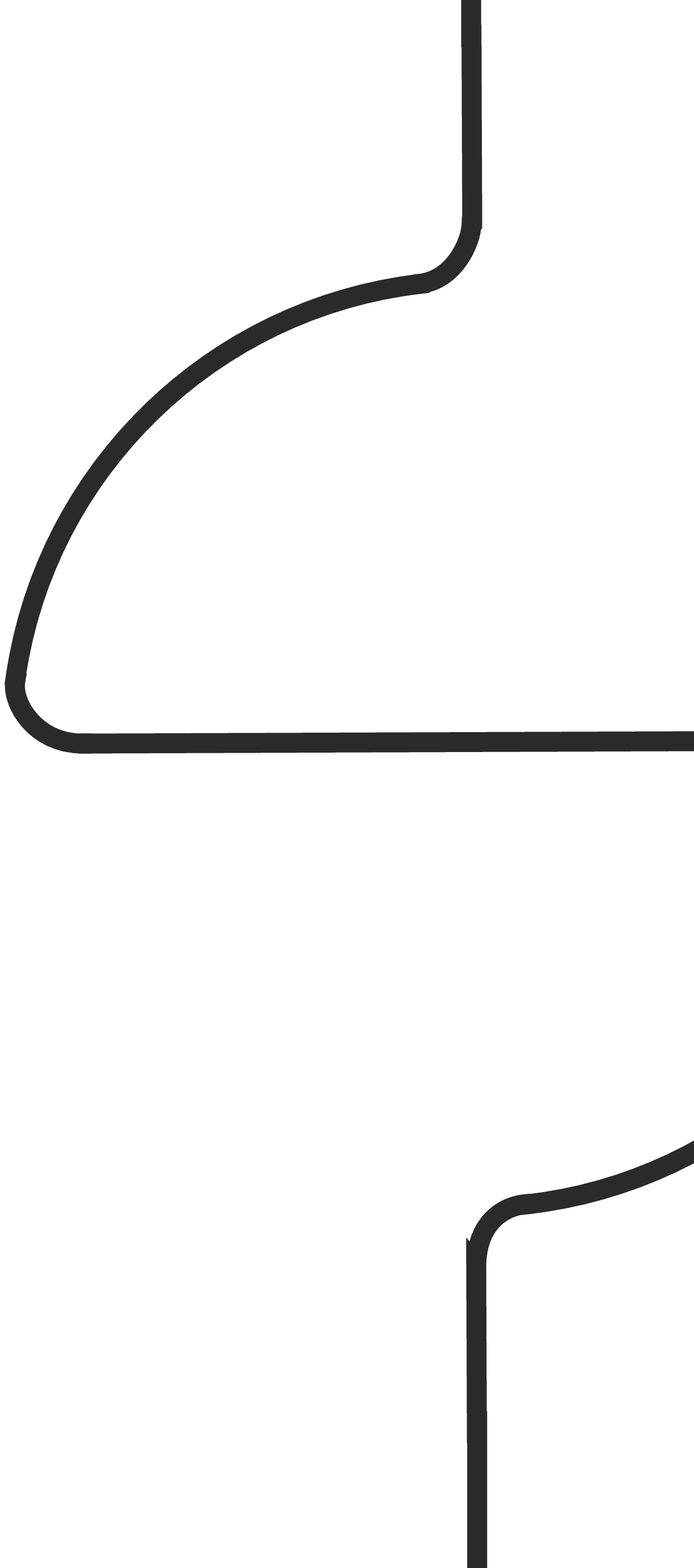}
\,\,\,\raisebox{1.5cm}{$\Longrightarrow$}
\,\,\,
 \includegraphics[height=2.8cm]{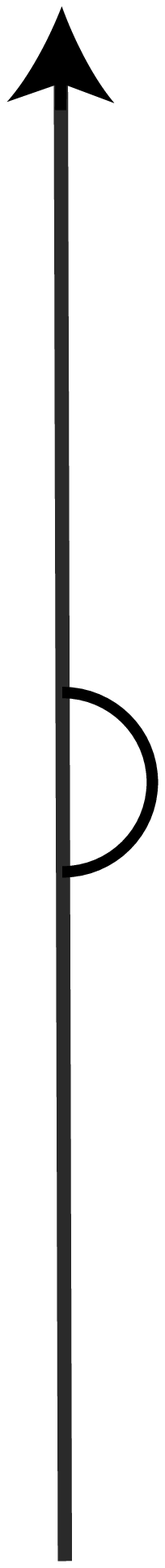} 
\end{center}
\caption{
The marker of the insertion of a leftward (resp. rightward) 
twisting arc}
\label{mark}
\end{figure}

From now on, we would like to twist a given planar 
curve by replacing a sufficiently small subarc 
with the above two twisting arcs.
Namely, one can attach the leftward (resp. rightward)
twisting arc into a given planar curve, and get a space curve. 
For the sake of simplicity, we indicate these two surgeries 
constructing space curves from a given planar curve 
symbolically as in Figure \ref{mark} left (resp. right).

If we connect two end points of a twisting arc by 
a planar arc in $xy$-plane, we get a closed curve.
Since the curvature function (as a plane curve) 
of a twisting arc near the two end points as a plane curve takes opposite
sign, the resulting closed curve has at least one inflection point.
We need such an operation to construct several M\"obius developables
in later. The existence of inflection points is really
needed for constructing rectifying M\"obius strips.
The following assertion is useful for counting M\"obius 
twisting number of our latter examples:

\begin{proposition}\label{prop:2-1}
Let $\tilde \sigma_+(t)$ and $\tilde \sigma_-(t)$ $(a\le t\le b)$ be
the leftward and rightward twisting arcs parameterizing the
set $\pi_+^{-1}(\hat \sigma)$
$($resp. $\pi_-^{-1}(\hat \sigma))$ respectively.
Then the space curves $\tilde \sigma_\pm(t)$
have no inflection points. 
Moreover, it holds that 
\begin{align}
&\op{Tw}_{\tilde \sigma_+}(D_+^\perp)
-\op{Tw}_{\tilde \sigma_+}(\mb e_3^\perp)
=\pi, \qquad
\op{Tw}_{\tilde \sigma_-}(D_-^\perp)
-\op{Tw}_{\tilde \sigma_-}(\mb e_3^\perp)
=-\pi, \label{eq:prop1}\\
&
\op{Tw}_{\tilde \sigma_+}(\eta_+)
-\op{Tw}_{\tilde \sigma_+}(\mb e_3^\perp)
=\pi,
 \qquad
\op{Tw}_{\tilde \sigma_-}(\eta_-)
-\op{Tw}_{\tilde \sigma_-}(\mb e_3^\perp)
=-\pi, \label{eq:prop2}
\end{align}
where $D_{\pm}(t)$ is the Darboux vector field of $\tilde \sigma_\pm(t)$,
$\mb e_3={}^t\!(0,0,1)$ and 
$$
\eta_\pm(t):=\tilde \sigma_\pm(t)\times \dot {\tilde \sigma}_\pm(t)
$$
is the (leftward) unit co-normal vector of $\tilde \sigma$
on the unit sphere $S^2$.
$($Here the normal sections 
 $D_{\pm}^\perp,\mb e_3^\perp$ 
with respect to ${\tilde\sigma}_\pm$
are obtained as the normal
parts of the vectors $D_{\pm},\mb e_3$.
See \eqref{eq:index}. $)$ 
\end{proposition}
\begin{proof}
It is sufficient to prove the
case of leftward twisting arc.
Let $\mb b(t)$ be the bi-normal vector of $\tilde \sigma_+$
as a space curve. Since $\tilde \sigma_+$ is a curve on
the unit sphere, the principal normal direction 
$\mb n(t)$ must be $-\sigma_+(t)$, 
and thus
$$
\mb b(t)=\mb t(t)\times \mb n(t)
=\tilde \sigma_+(t)\times \mb t(t)=\eta_+(t),
$$
where $\mb t(t):=\dot\gamma(t)/|\dot\gamma(t)|$.
Moreover, by the definition of the normalized 
Darboux vector field $D_+(t)$, we have
$$
D^\perp_+(t)=\mb b(t)=\eta_+(t).
$$
Thus the first formula reduces to the
second one.\\
Let $\theta(t)$ be the smooth function
which gives the leftward angle of 
$\eta_+(t)$ from $\mb e_3^\perp$.
Then, we have
$$
\op{Tw}_{\tilde \sigma_+}(\eta_+)
-\op{Tw}_{\tilde \sigma_+}(\mb e_3^\perp)
=\theta(1)-\theta(0).
$$
Let $\mb t(t)$ be the unit tangent vector
of $\tilde \sigma_+$ as a space curve.
Then by definition of $\tilde \sigma_+$,
we have 
$$
\mb t(0)=\mb t(1),\qquad \mb n(0)=-\mb n(1)
$$
which yield
\begin{equation}\label{eq:b1}
\eta_+(0)=\mb b(0)=\mb t(0)\times \mb n(0)=
-\mb t(1)\times \mb n(1)=-\mb b(1)=-\eta_+(1).
\end{equation}
On the other hands, $\tilde \sigma_+(t)$
lies in $xy$-plane near $t=a,b$, the vector 
$\eta_+(t)=\mb b(t)$ is proportional
to $\mb e_3$ there.
Thus we have that 
$$
\theta(1)-\theta(0)=\pi \,\, \mod 2\pi \Z.
$$
Since we can easily check that $\theta(t)\ge 0$,
we get $\theta(1)-\theta(0)=\pi$,
which proves \eqref{eq:prop2}.
\end{proof}

\begin{lemma}\label{lem2.2}
Let $\gamma(t)$ be a spherical curve parametrized 
by the arclength parameter.
Then the leftward conormal vector field
$$
\eta(t):=\gamma(t)\times \dot\gamma(t)
$$
is parallel with respect the normal connection
of $\gamma(t)$. In particular,
$
F_{\gamma,\eta}(t,u) 
$
is a principal developable strip.
\end{lemma}

\begin{proof}
A normal vector field $\xi(t)$
along $\gamma$ is parallel with respect to 
the normal connection if and only if
$\dot \xi(t)$ is proportional to $\dot \gamma$.
Applying the Frenet formula, we have
$$
\dot\eta(t)=\gamma(t)\times \ddot\gamma(t)
=\kappa(t) \gamma(t)\times \mb n(t),
$$
where $\mb n(t)$ and $\kappa(t)$ 
is the principal normal vector and
the curvature function of $\gamma(t)$ as a
space curve.
Since $\gamma$ and $\mb n$ are both perpendicular to
$\dot \gamma$, the vector $\gamma\times \mb n$
is proportional to $\dot \gamma$,
which proves the assertion.
\end{proof}

\begin{definition}\label{def:rect}
Let $\tilde \sigma_+(t)$(resp. $\tilde \sigma_-(t)$) be
the leftward  (resp. rightward) twisting arc as in Proposition 
\ref{prop:2-1}.
Then 
$$
F^{\pm}_{p}(t,u):=\tilde \sigma_\pm(t)+u\eta_{\pm}(t)
\qquad (\eta_\pm(t):=\tilde\sigma_\pm(t)\times \dot{\tilde\sigma}_\pm(t))
$$
is called the {\it principal twisting strip}
and
$$
F^\pm_{g}(t,u):=\tilde \sigma_\pm(t)+u D_{\pm}(t)
$$
is called the {\it rectifying twisting strip}, where
$D_{\pm}(t)$ is the normalized Darboux field of $\tilde \sigma_{\pm}$.
\end{definition}

By Proposition \ref{prop:2-1} and Lemma \ref{lem2.2},
$F^\pm_{p}$ is a 
principal developable satisfying \eqref{eq:prop2},
and $F^\pm_{g}$ is a rectifying developable
satisfying \eqref{eq:prop1}.

\begin{figure}[htb]
\begin{center}
 \includegraphics[width=2.0cm]{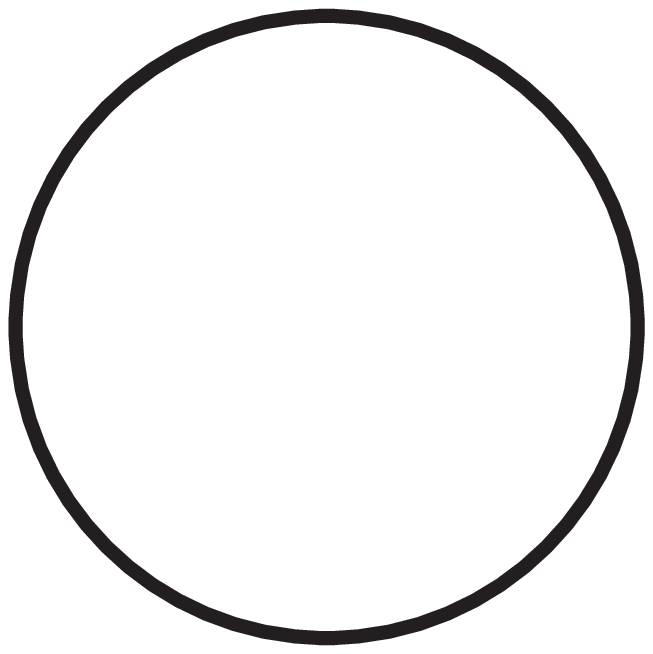}
\,\,\,\raisebox{1.0cm}{$\Longrightarrow$}
\,\,\,
 \includegraphics[width=2.2cm]{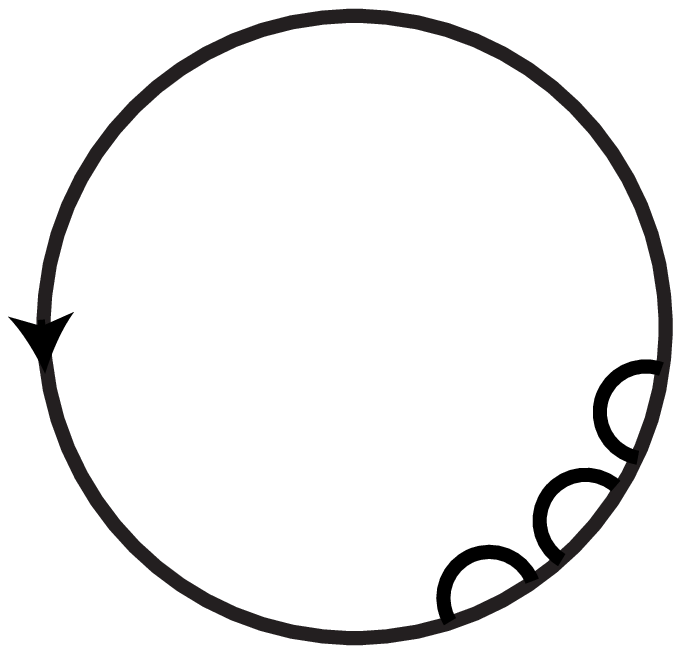} 
\,\, \,
 \includegraphics[width=2.2cm]{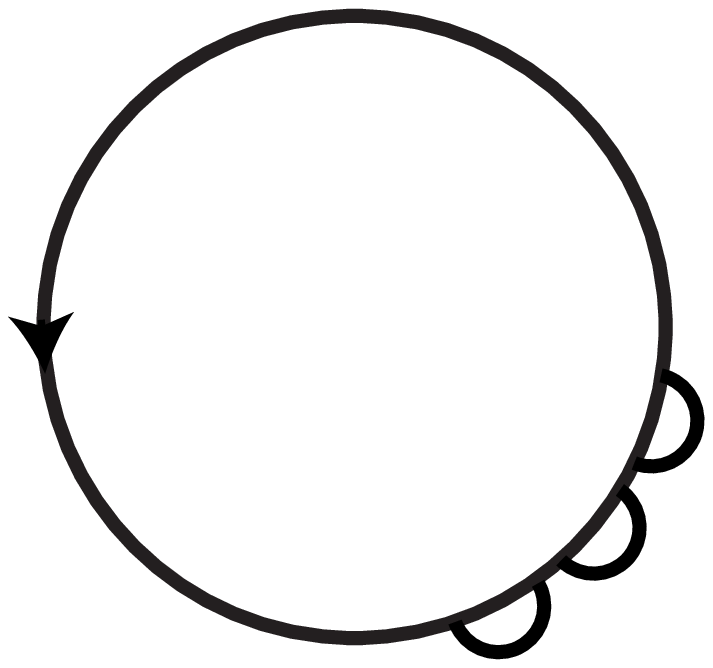} 
\end{center}
\caption{The construction of $C_{2m+1}$
via $C$.}
\end{figure}

\begin{theorem}\label{thm:smooth}
For an arbitrarily given isotopy type of 
M\"obius strip, there exists a $C^\infty$ 
principal $($resp. rectifying$)$
M\"obius developable in the same isotopy class.
\end{theorem}

\begin{proof}
First, we construct an unknotted principal M\"obius developable
of a given M\"obius twisting number from a circle:
Consider a circle $C$ in the $xy$-plane.
We insert $2m+1$ leftward (resp. rightward) twisting arcs
into $C$ and denote it by $C_{2m+1}$ or $C_{-2m-1}$
(See Figure 4.). 
If we build $2m+1$ principal twisting strips
(each of which is congruent to $F^\pm_p$)
on these twisting arcs,
then we get a principal $C^\infty$ M\"obius developable
$F_{2m+1}$ whose centerline is $C_{2m+1}$.
(Let $\gamma(t)$ ($a\le t \le b$) be a parametrization of
centerline of $F_{2m+1}$. Then we can write
$$
F_{2m+1}(t,u)=\gamma(t)+u P(t)\qquad (a\le t \le b,\,\,|u|<\epsilon).
$$
The image of the center line $\gamma(t)$ is a union of 
$m$ planar arcs and $m$ twisting arcs.
On each planar arcs $P(t)$
is equal to $ \mb e_3={}^t\!(0,0,1)$.
On the other hand, $P(t)$ coincides with the co-normal
vector on each twisting arc as a spherical curve. 
Since the twisting arc is planar near two end points,
$P(t)$ is smooth at each end points of
twisting arcs. Consequently,
$P(t)$ satisfies the condition of 
Proposition \ref{prop:principal} such that
$P(a)=-P(b)$.)
By \eqref{eq:index} and \eqref{eq:prop2}, 
the M\"obius twisting number of $F_{2m+1}$ 
is equal to $-(2m+1)/2$ (resp. $(2m+1)/2$)
if we insert the leftward (resp. rightward)
twisting strips.

\medskip
Instead of principal twisting strips, we can insert
rectifying twisting strips $F^\pm_g$ into $C_{2m+1}$.
Then by \eqref{eq:index} and \eqref{eq:prop1},
we also get a rectifying $C^\infty$
M\"obius developable 
with the M\"obius twisting number $\pm (2m+1)/2$.

\medskip
Next, we construct a knotted principal M\"obius developable
of a given M\"obius twisting number via a knot diagram.
It should be remarked that
the isotopy type of the given embedded M\"obius strip
is determined by its M\"obius twisting number and the
knot type of its centerline. (See \cite{RR2}.)
Let $\gamma$ be the planar curve corresponding to
the diagram.
We replace every crossing of $\gamma$ by
a pair of leftward and rightward twisting arcs
as in Figure \ref{crossing1} (right).
For the sake of simplicity, we indicate this operation
as in Figure \ref{crossing2}.
When we will accomplish to construct the associated 
M\"obius developable, this operation as in Figure 5
does not effect the M\"obius 
twisting number, since
the signs of the two twisting arcs are opposite.

\begin{figure}[htb]
\begin{center}
 \includegraphics[height=3.2cm]{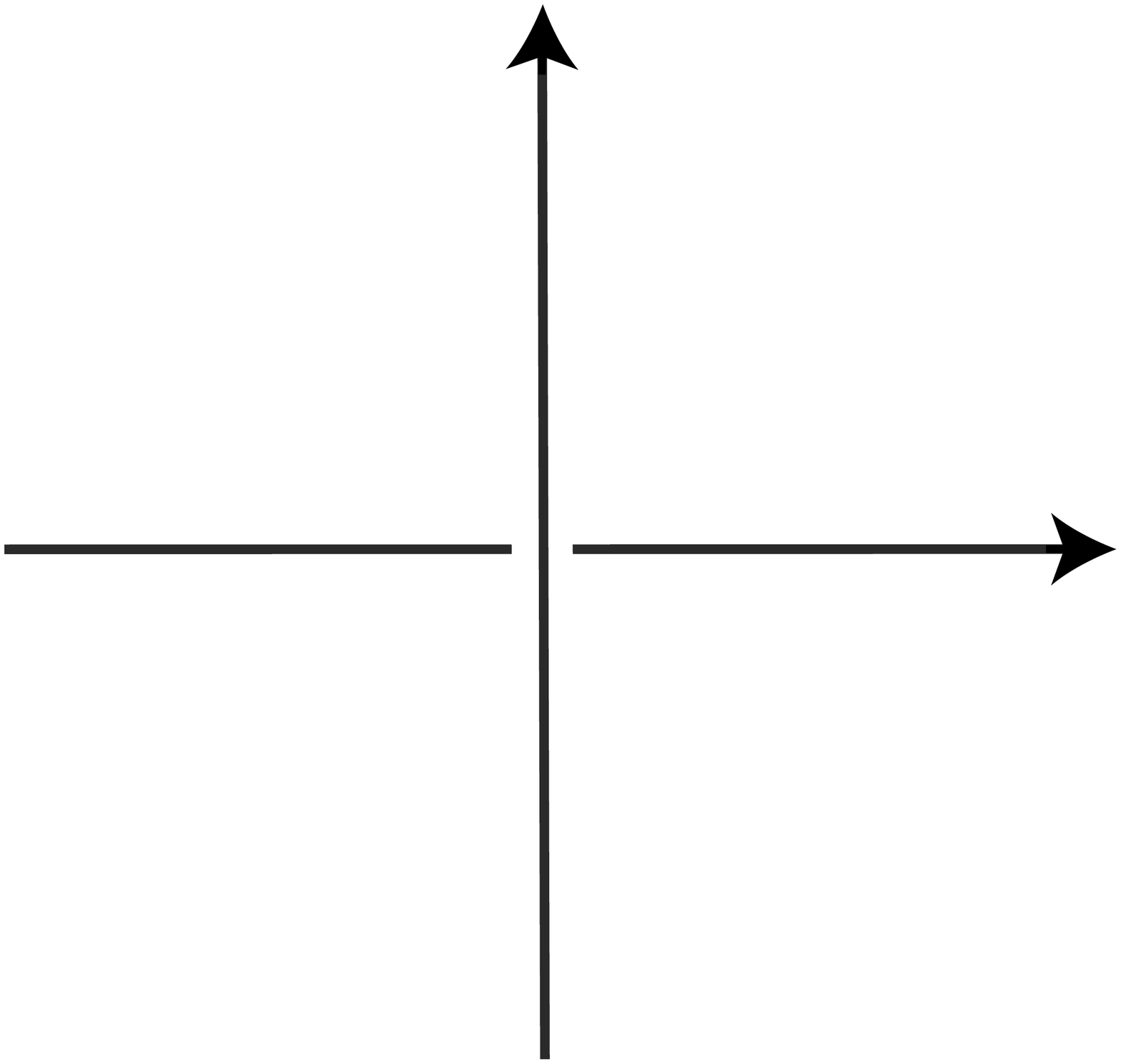}
\,\,\,\raisebox{1.5cm}{$\Longrightarrow$}
\,\,\,
 \includegraphics[height=3.5cm]{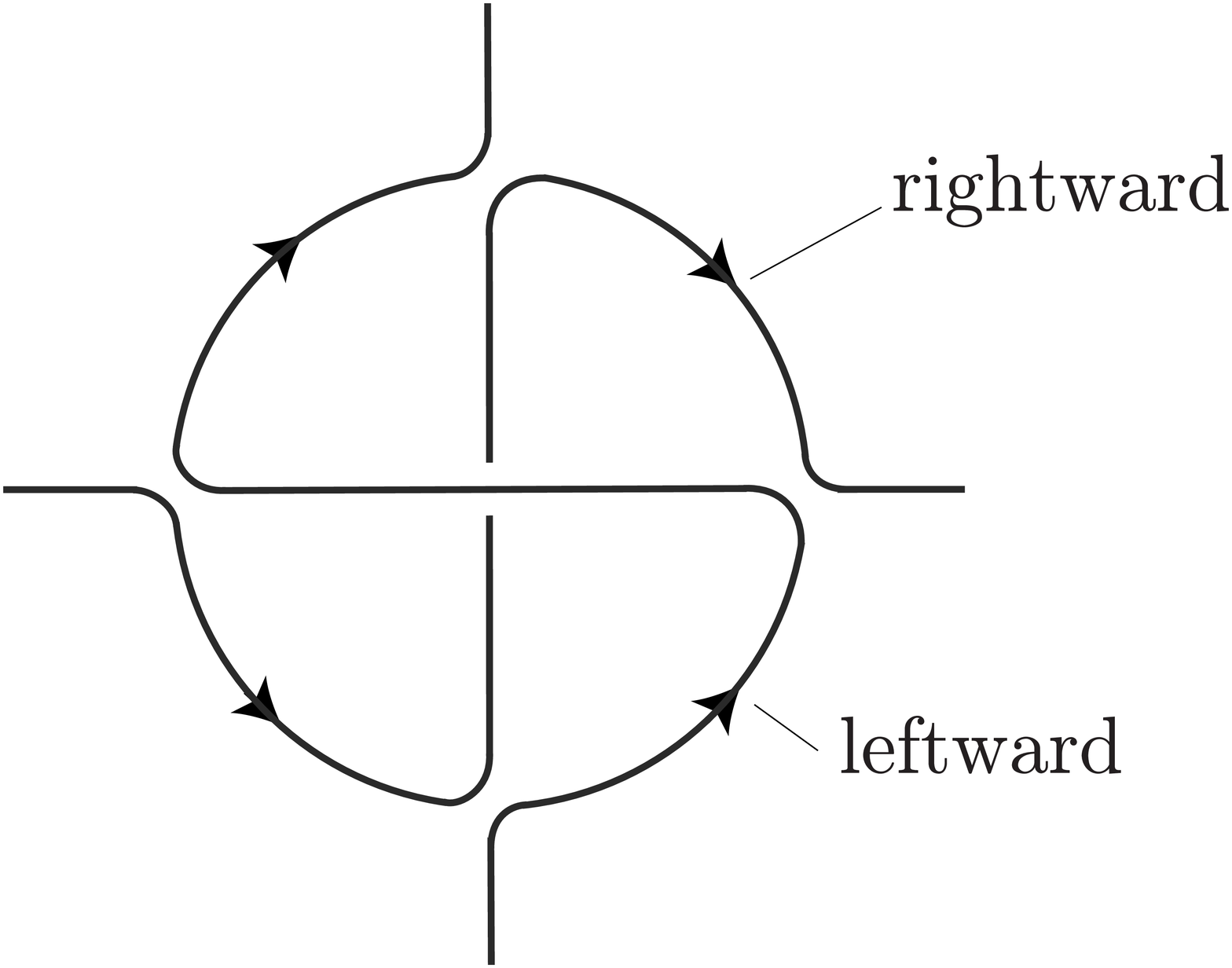} 
\end{center}
\caption{The crossing with the pair of twisting arcs}
\label{crossing1}
\end{figure}

\begin{figure}[htb]
\begin{center}
 \includegraphics[height=2.5cm]{crossing2.eps}
\,\,\,\raisebox{1.0cm}{$\Longrightarrow$}
\,\,\,
 \includegraphics[height=2.5cm]{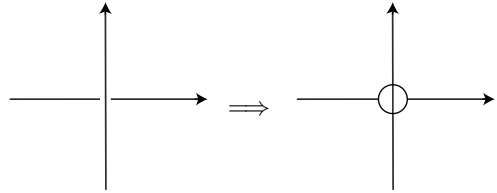} 
\end{center}
\caption{The marker of the pair of twisting arcs at a crossing 
}
\label{crossing2}
\end{figure}

For example, letting $K$ be a knot 
diagram $3_1$ of the trefoil knot
as in Figure \ref{trefoil} left,
we replace each crossing by
 a pair of leftward and rightward twisting 
arcs (as in Figure 5 and Figure 6), and insert $2m+1$  
leftward (resp. rightward) 
twisting arcs as in Figure \ref{trefoil} right.
Then we get an embedded closed space curve $C^{K}_{2m+1}$
($m\in \Z$) which is isotopic to the knot $K$.
If we build  principal twisting strips
on all of the twisting arcs we inserted,
then we get a principal $C^\infty$ M\"obius 
developable $F^K_{2m+1}$.
Since all crossing of $3_1$ are positive,
the writhe is $3$, and thus
the formula \eqref{eq:index} and \eqref{eq:prop2} yields that
the M\"obius twisting number of $F^K_{2m+1}$ is $3\mp (2m+1)/2$.
Since $m$ is an arbitrary non-negative integer, we prove the 
existence of principal  M\"obius strip
for the case of trefoil knot.

Similarly,
we can prove the
existence of principal  M\"obius strip
$F^K_{2m+1}$ with an arbitrary M\"obius twisting number
for an arbitrary given
knot diagram $K$.

\begin{figure}[htb]
\begin{center}
 \includegraphics[height=3.0cm]{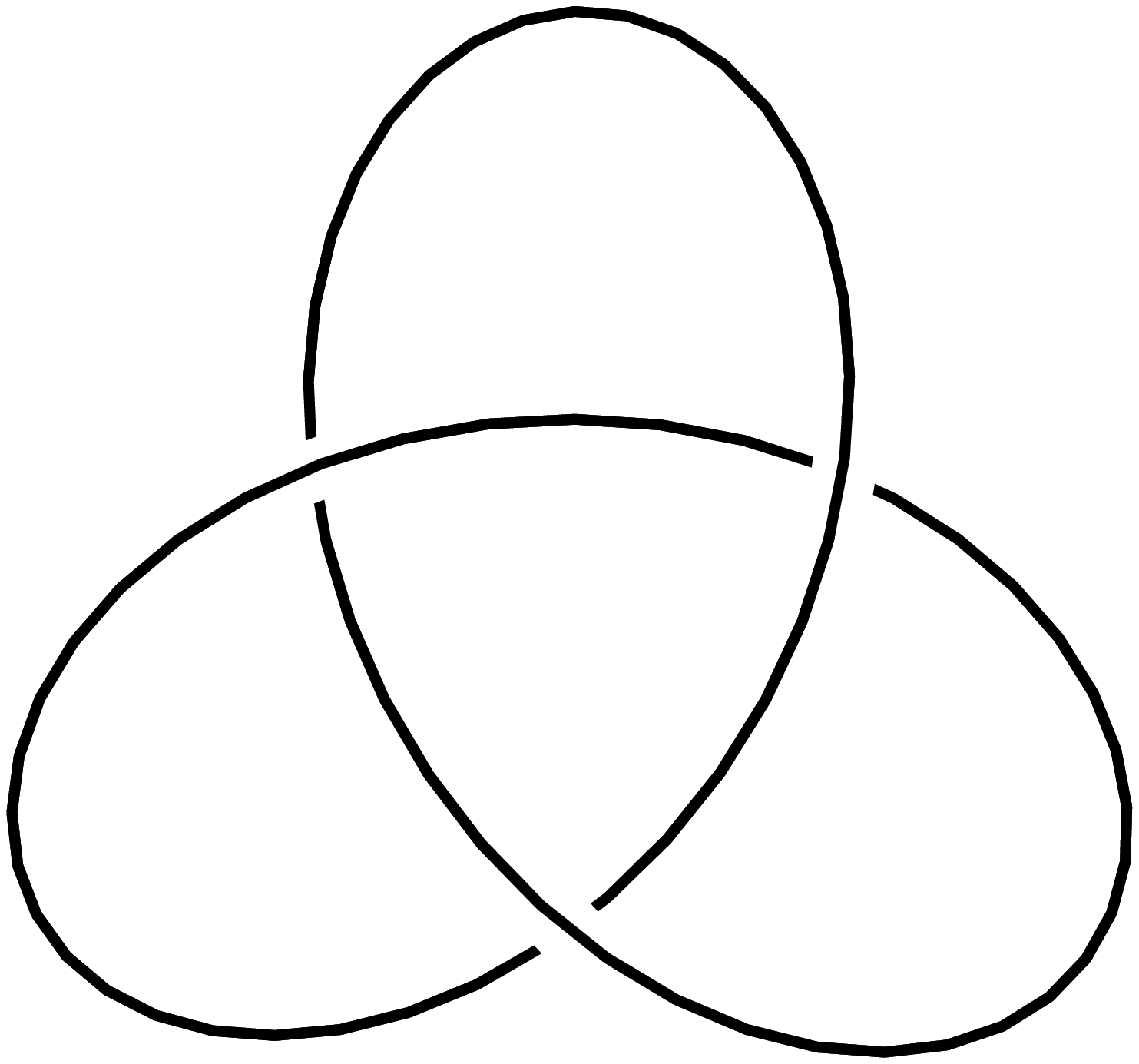}
\,\,\,\raisebox{1.5cm}{$\Longrightarrow$}
\,\,\,
 \includegraphics[height=3.2cm]{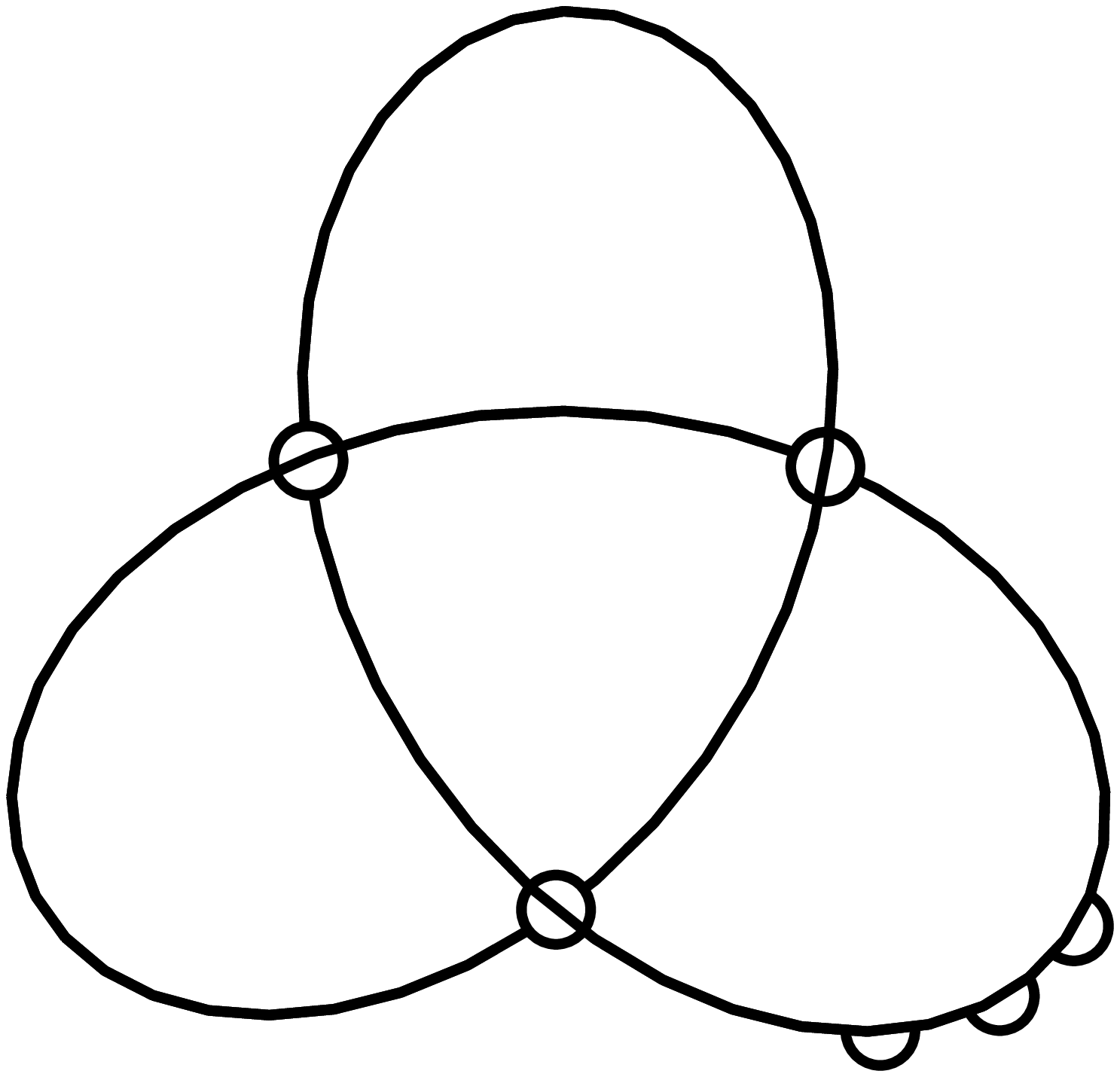} 
\end{center}
\caption{The construction of a M\"obius developable via the 
knot diagram of $3_1$}
\label{trefoil}
\end{figure}

\medskip
Instead of the principal twisting strips, we can
insert the  rectifying twisting strips (cf. Definition
\ref{def:rect}).
Then we also get a rectifying $C^\infty$ M\"obius developable
with an arbitrary isotopy type at the same time.
\end{proof}

\medskip
\noindent
{\bf (Properties of 
asymptotic completion of M\"obius strips)}\\
Let $M^2$ be a $2$-manifold and 
$f:M^2\to \mathbf R^3$ a $C^\infty$-map.
A point $p\in M^2$ is called {\it regular} if
$f$ is an immersion on a sufficiently small neighborhood
of $p$, and is called {\it singular} if it is not regular. 
Moreover, $f:M^2\to \mathbf R^3$ 
is called a {\it (wave) front} if
\begin{enumerate}
\item there exists a unit vector field 
$\nu$ along $f$ such that $\nu$
is perpendicular to the image of
tangent spaces $f_*(TM)$.
( $\nu$ is called the
{\it unit normal vector field} of $f$, which 
can be identified with the Gauss map
$\nu:M^2\to \mathbf R^3$. )
\item 
The pair of maps
$$
L:=(f,\nu):M^2\to \mathbf R^3\times S^2 (\cong T^*_1\mathbf R^3)
$$
gives an immersion.
\end{enumerate}
On the other hand, a smooth map $f:M^2\to \mathbf R^3$ 
is called a {\it p-front} if it is locally a front, that is,
for each $q\in M^2$, there
exists an open neighborhood $U_q$ such that
the restriction $f|_{U_q}$ gives a front.
By definition, a front is a p-front.
A p-front is a front if and only if it has globally
defined unit normal vector fields (namely, it is
 co-orientable).
\begin{definition}\label{def:complete}
(\cite{MU})\,
The first fundamental form 
$ds^2$ of 
a flat p-front $f:M^2\to \R^3$  
is called {\it complete} if
there exists a symmetric covariant tensor $T$ on $M^2$
with compact support
such that $ds^2+T$ gives a complete metric
on $M^2$. 
On the other hand, $f$ is called {\it weakly complete}
if the sum of the first fundamental form and the
third fundamental form 
$$
ds^2_\#:=df\cdot df+d\nu\cdot d\nu
$$
gives a complete Riemannian metric on $M^2$.
\end{definition}
A front is called {\it flat} if $\nu:M^2\to S^2$ is
degenerate everywhere. 
Parallel surfaces $f_t\, (t\in \R)$ 
and the caustic $\mathcal C_f$ of a
flat front $f$ are all flat. 
Weakly complete flat p-front is complete if and
only it is weakly complete and the singular set is compact.
(See \cite[Corollary 4.8]{MU}.)
Let $\epsilon>0$ and
$$
F(=F_{\gamma,\xi}(t,u))=\gamma(t)+u \xi(t) \qquad (|u|<\epsilon),
$$
be a flat M\"obius developable
defined on a closed interval $t\in [a,b]$.
Then
$$
\tilde F(t,u)=\gamma(t)+u \xi(t) \qquad (u\in \R)
$$
as a map of $S^1\times \R$ is called
the {\it asymptotic completion} of $f$.
We can prove the following:
\begin{corollary}\label{cor:complete}
For an arbitrary given isotopy type of M\"obius
strip, there exists a principal M\"obius developable $f$
in the same isotopy class
whose asymptotic completion $\tilde f$ gives a weakly complete
flat p-front. 
\end{corollary}
In \cite[Theorem A]{MU}, it is shown that complete flat p-front
is orientable. In particular, the singular set of $\tilde f$
as above cannot be compact.
\begin{proof}
Let $F$ be a principal M\"obius strip constructed in
the proof of Theorem \ref{thm:smooth}.
We can write
$$
\tilde F(t,u)=\gamma(t)+u P(t) \qquad (t\in [a,b],\,\, u\in \R),
$$
where $\gamma(t)$ be the embedded space curve $C_{2m+1}$ or
$C^K_{2m+1}$.
By taking $t$ to be the arclength parameter of $\gamma$,
we may assume 
\begin{equation}\label{eq:unit1}
|\dot \gamma(t)|=1\qquad (t\in [a,b]).
\end{equation}
Since $F$ is principal, the asymptotic direction
$P(t)$ is parallel with
respect to normal section. 
In particular, we may also assume that
\begin{equation}\label{eq:unit2}
|P(t)|=1\qquad (t\in [a,b]),
\end{equation}
and
\begin{equation}\label{eq:unit3}
\dot P(t)=\lambda(t)\dot \gamma(t) \qquad (t\in [a,b]).
\end{equation}
As seen in the proof of Theorem \ref{thm:smooth}, 
we may assume there exist points 
$$
a< p_1<  q_1 < p_2<  
q_2<  \cdots < p_n\le  q_n <b
$$
such that the interval $(p_j,q_j)$ corresponding to
the twisting arcs, in particular, we have 
\begin{enumerate}
\item The open subarc $\gamma(t)$ 
$(t\in \dy\bigcup_{j=1}^n (p_j,q_j))$ has no 
inflection points as a space curve,
\item $P(t)=\mb e_3$ for $t\not\in \dy\bigcup_{j=1}^n (p_j,q_j)$.
\end{enumerate}
As seen in the proof of theorem \ref{thm:smooth},
the curve $\gamma$ is constructed from a knot diagram $K$.
We set  
$$
\nu(t):=\dot \gamma(t)\times P(t).
$$
Then it gives the normal vector of $F(t,u)$.
If we choose the initial knot diagram generically, we may
assume that the number of inflection points
on the diagram is finite. Then we can insert principal twisting arcs
in the diagram apart from these inflection points.
Since $\gamma$ is principal, the Weingarten formula yields that
$\dot \nu(t)$ gives a principal direction (cf. \eqref{eq:new}), and
$|\dot \nu(t)|$ gives the absolute value of
the principal curvature function of $f$.
So $|\dot \nu(t)|$ does not vanish if $t$ 
is not an inflection point of $\gamma(t)$.
Thus there exists a positive constant $\rho_0(<1)$ such that
$$
|\dot\nu(t)|\ge \rho_0 \qquad (t\in \dy\bigcup_{j=1}^n (p_j,q_j)).
$$
Since $P(t)$ is perpendicular to $\dot \gamma(t)$,
\eqref{eq:unit1}, \eqref{eq:unit2} and \eqref{eq:unit3}
yields that
$$
ds^2_\#=ds^2+d\nu^2=\biggl((1+u\lambda(t))^2dt^2+du^2\biggr)+|\dot\nu(t)|^2dt^2.
$$
Then we have that
\begin{equation}\label{eq:metric1}
ds^2_\#\ge du^2+|\dot\nu(t)|^2dt^2\ge du^2+|\rho_0|^2dt^2
\qquad (t\in \dy\bigcup_{j=1}^n (p_j,q_j)).
\end{equation}
Next we suppose that $t\not \in \dy\bigcup_{j=1}^n (p_j,q_j)$.
Then $P(t)=\mb e_3$ holds and thus
$\lambda(t)$ vanishes.
Since $\rho_0<1$, we have
\begin{equation}\label{eq:metric2}
ds^2_\#=(dt^2+du^2)+|\dot\nu(t)|^2dt^2\ge (dt^2+du^2)\ge 
du^2+|\rho_0|^2dt^2.
\end{equation}
By \eqref{eq:metric1} and \eqref{eq:metric2},
we have $ds^2_\#\ge du^2+|\rho_0|^2dt^2$ for all 
$t\in [a,b]$. In particular, $ds^2_\#$ is positive definite
and $\tilde f$ is a front.
Moreover, since $du^2+|\rho_0|^2dt^2$ is a complete
Riemannian metric on $S^1\times \R$, so is $ds^2_\#$,
which proves the assertion.
\end{proof}
\medskip
\noindent
({\it Proof of Theorem A}.)\\
Let $F$ be a principal M\"obius strip constructed as in
the proof of Corollary \ref{cor:complete}, that is
we can write
$$
F(t,u)=\gamma(t)+u P(t) \qquad (t\in [a,b],\,\, |u|<\epsilon).
$$
We fix an integer $m\in \Z$ arbitrarily.
Then we can take $F$ so that
\begin{equation}\label{eq:TwP}
\op{Tw}_\gamma(P)=\frac{2m+3}2.
\end{equation}
Moreover, we may assume that 
$$
a=0,\qquad b=2\pi.
$$ 
Here  $\gamma$ lies on $xy$-plane
when $t\not \in \dy\bigcup_{j=1}^n (p_j,q_j)$.
So without loss of generality,
we may also assume that $0\not \in \dy\bigcup_{j=1}^n (p_j,q_j)$.
Then $P(t)$ is uniquely determined by the 
initial condition $P(0)=\mathbf e_3$.
Let
$$
\Pi:\R^3\to \R^2
$$
be the projection into $xy$-plane.
We set
$$
\gamma_d(t):=(1-d)\gamma(t)+d 
\,\Pi\circ \gamma (t)
\qquad (0\le d\le 1).
$$
Then $\gamma_d$ has same isotopy type
as $\gamma=\gamma_0$ for each $d\in (0,1)$.
Consider the Fourier expansion of $\gamma_d(t)$
under the identification $S^1=\R/(2\pi \Z)$
$$
\gamma_d(t)=a_0(d)+\sum_{n=1}^\infty
\biggl( 
a_n(d) \cos(nt) + b_n(d) \sin(nt)
\biggr),
$$
and let
$$\gamma_{d,n}(t)
:=a_0(d)+\sum_{j=1}^n
\biggl( 
a_j(d) \cos(jt) + b_j(d) \sin(jt)
\biggr)
\qquad (n=1,2,3,...)
$$
be the $n$ th approximation of $\gamma_d(t)$.
Then $\{\gamma_{d,n}\}$ is a
family real analytic curves
uniformly converges to $\gamma_d$.
Since $d$ is a real analytic parameter of $\gamma_d$,
$$
a_0(d),a_1(d),b_1(d),a_2(d),b_2(d),\cdots
$$
are all real analytic functions of $d$.
For each positive integer $n$ and $d\in [0,1]$, there exists
a unique vector field $P_{d,n}(t)$ along $\gamma$
such that $P(0)=\mb e_3$ and $\dot P(t)$ is proportional to
$\dot \gamma$.
Moreover,
$$
\lim_{n\to \infty}P_{0,n}(t)=P(t)
$$
and
$$
\lim_{n\to \infty}P_{1,n}(t)=\mb e_3.
$$
Since $\gamma_{1,n}$ is a plane curve in $xy$-plane,
we have
$$
\lim_{n\to \infty}\op{Tw}_{\gamma_{0,n}}(P_{0,n})=\frac{2m+3}2,
\qquad
\lim_{n\to \infty}\op{Tw}_{\gamma_{1,n}}(P_{1,n})=0.
$$
By the intermediate value theorem, there exists
$\delta_0\in (0,1)$ such that
$$
\op{Tw}_{\gamma_{\delta_0,n}}(P_{\delta_0,n})=\frac{2m+1}2,
$$
for sufficiently large $n$.
By \eqref{eq:index}, 
$$
F_n(t,u):=\gamma(t)+u P_{\delta_0,n}(t)
$$
gives a real analytic principal M\"obius strip
of twisting number 
$$
-\frac{2m+1}2+\op{Tw}_\gamma(\mb e_3^\perp)+\op{Wr}_{\mb e_3}(K)
$$
where $\op{Wr}_{\mb e_3}(K)$ is the writhe of the knot diagram
$K$. (If $K$ is un-knotted, the writhe vanishes.)
Since $\op{Tw}_\gamma(\mb e_3^\perp)$ and $\op{Wr}_{\mb e_3}(K)$
are both fixed integers and $m\in \Z$ is arbitrary,
this $F_n$ gives the desired real analytic principal M\"obius strip.
\hfill \qed

\section{Proof of Theorem B.}

We construct a real analytic M\"obius developable,
by a deformation of a $C^\infty$ M\"obius developable.
For this purpose, the rectifying $C^\infty$ M\'obius developables
given in the previous section is not sufficient 
and we prepare the following proposition instead:
(In fact, we must control inflection points on 
the centerline much more strictly to apply 
Corollary \ref{cor:generic}.)

\begin{proposition}\label{thm:admissible}
There exists a rectifying $C^\infty$ M\"obius developable
with an arbitrary isotopy type such that its centerline
$$
\gamma(t)=(x(t),y(t),z(t))\qquad (|t|\le \pi )
$$
as a $2\pi$-periodic embedded space curve satisfies
\begin{enumerate}
\item $\gamma(t)$ has a unique inflection point at $t=0$,
namely, $\dot \gamma(t)\times \ddot\gamma(t)\ne 0$ holds
for $t\ne 0$,
\item $\dot y(0)=\ddot y(0)=0$ and $\dddot y(0)\ne 0$,
\item $\dot z(0)=\ddot z(0)=\dddot z(0)=z^{(4)}(0)=0$.
\end{enumerate}
In particular, $t=0$ is the generic inflection point 
such that $\op{det}(\dot\gamma(0),\gamma^{(3)}(0),
\gamma^{(4)}(0))=0$ $($cf. Corollary \ref{cor:generic}$)$.
\end{proposition}

\begin{figure}[htb]
\begin{center}
 \includegraphics[height=2.0cm]{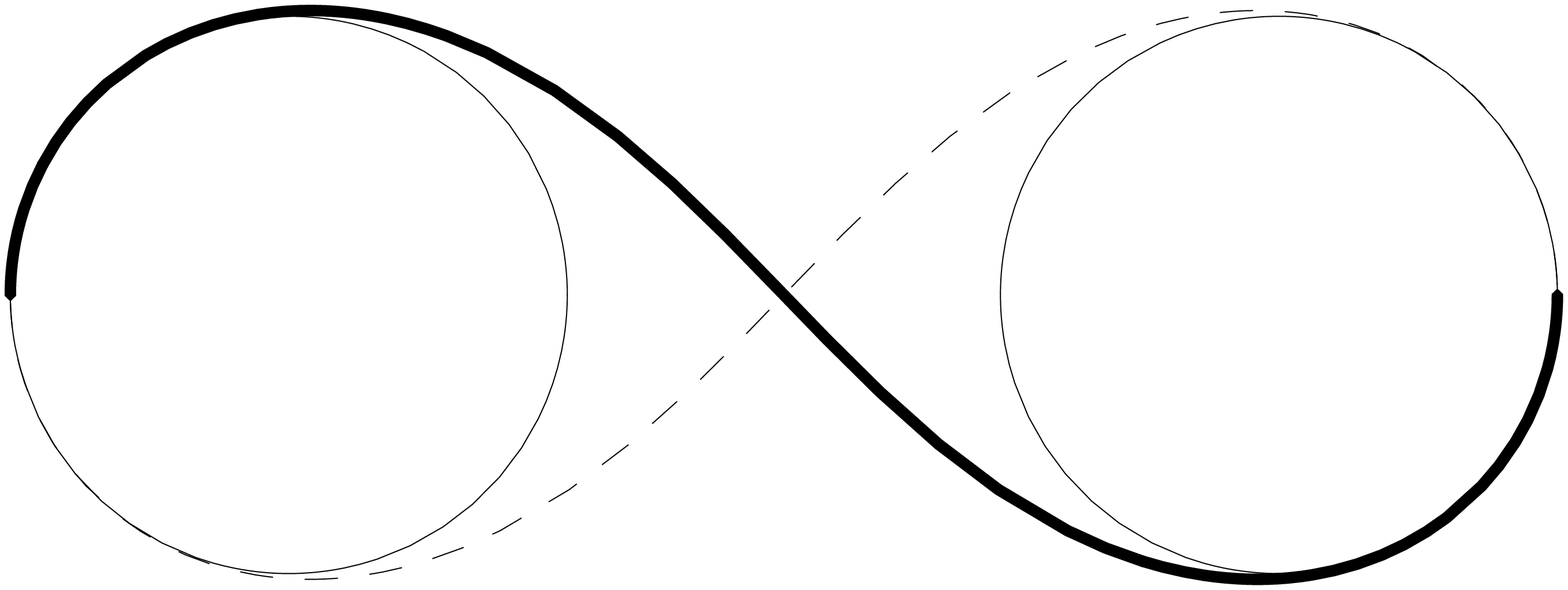}
\qquad
 \includegraphics[height=2.0cm]{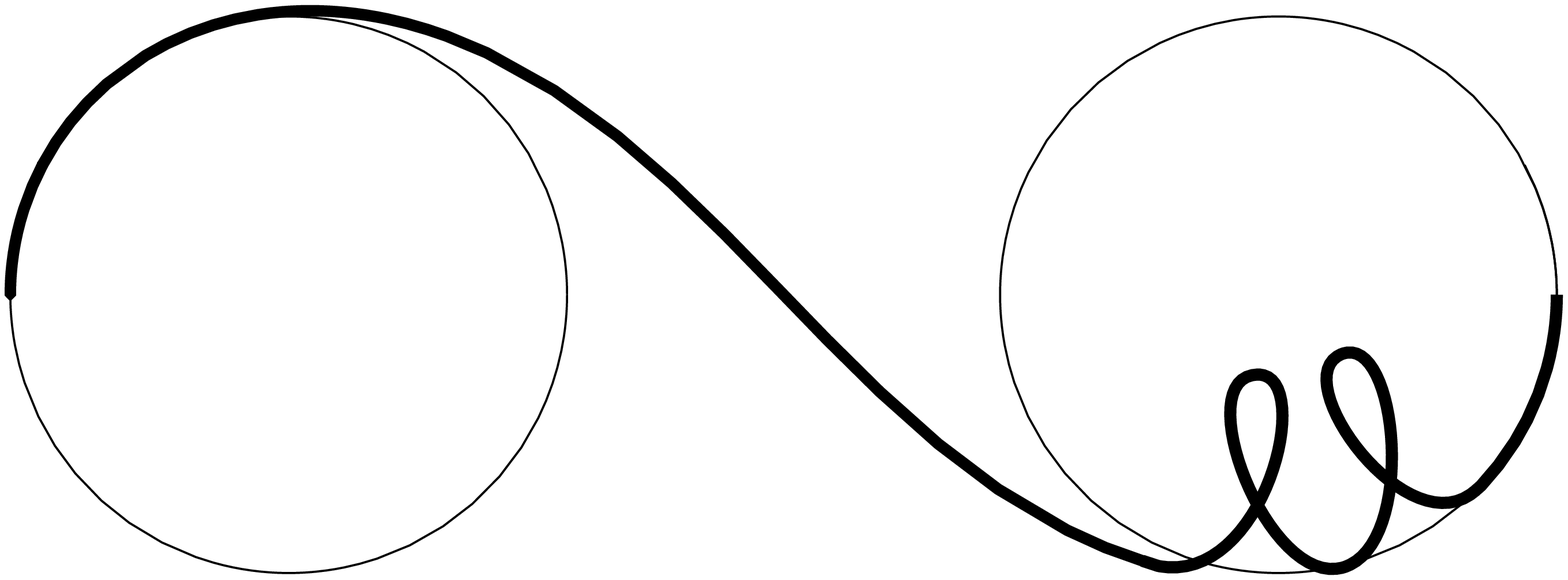}
\end{center}
\caption{$\sigma_0$ and \lq $\sigma_0$ with two loops \rq }
\label{cate23}
\end{figure}

To prove the proposition, we need additional 
special arcs in $\R^3$:


\medskip
\noindent
{\bf (The $S$-arc)}
The map  
$$
t\mapsto \frac{\cos t }{1 + \sin^2 t}\pmt{1\\
\sin t}
\qquad (0 \le t\le 2\pi)
$$
parametrizes a lemniscate given by
$$
(x^2 + y^2)^2 = x^2 - y^2
$$
as in Figure \ref{cate23} (left)
in the $xy$-plane.
The osculating conics at $t=0,\pi$ are 
exactly two circles
$$
(x\pm a)^2+y^2=b^2,
$$
which are inscribed in the lemniscate 
and meet the lemniscate with $C^3$-regularity, where
\begin{equation}\label{eq:ab}
a=\frac{2}{3},\qquad b=\frac{1}{3}.
\end{equation} 
So we set
$$
\gamma(t):=\frac{\cos t }{1 + \sin^2 t}
\pmt{1\\ \sin t}
\qquad (\pi \le t\le 2\pi).
$$
Since $\gamma(t)$ has $C^3$-contact with the osculating
circles $C_{\pi}$ and $C_{2\pi}$ at $t=\pi,2\pi$, 
we can give a $C^3$-differentiable perturbation
of  $\gamma$ near $t=\pi,2\pi$ such that
the new curve $\sigma_0(t)$ ($\pi\le t\le 2\pi$)
after the operation has $C^\infty$-contact with the 
circles $C_{\pi}$ and $C_{2\pi}$.
This new curve $\sigma_0$ 
is called the {\it $S$-arc} as in Figure \ref{cate23}
(left).

\begin{figure}[htb]
\begin{center}
 \includegraphics[height=2.0cm]{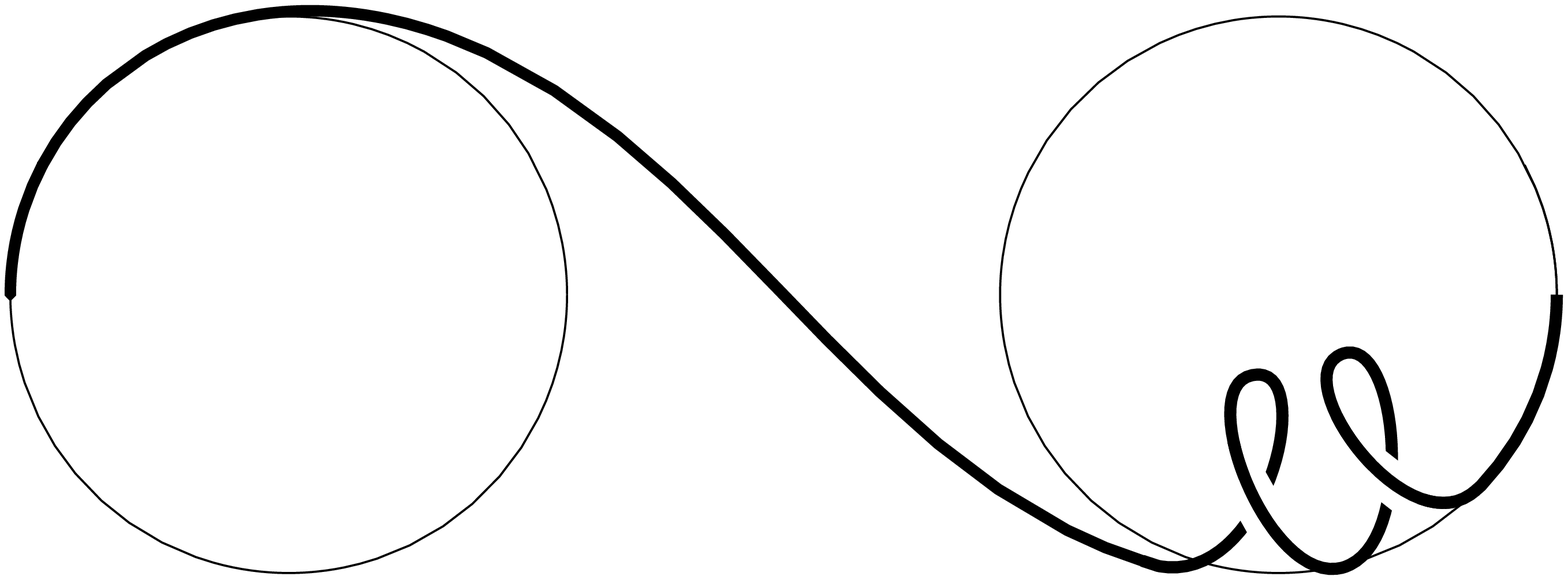}
\qquad
 \includegraphics[height=2.0cm]{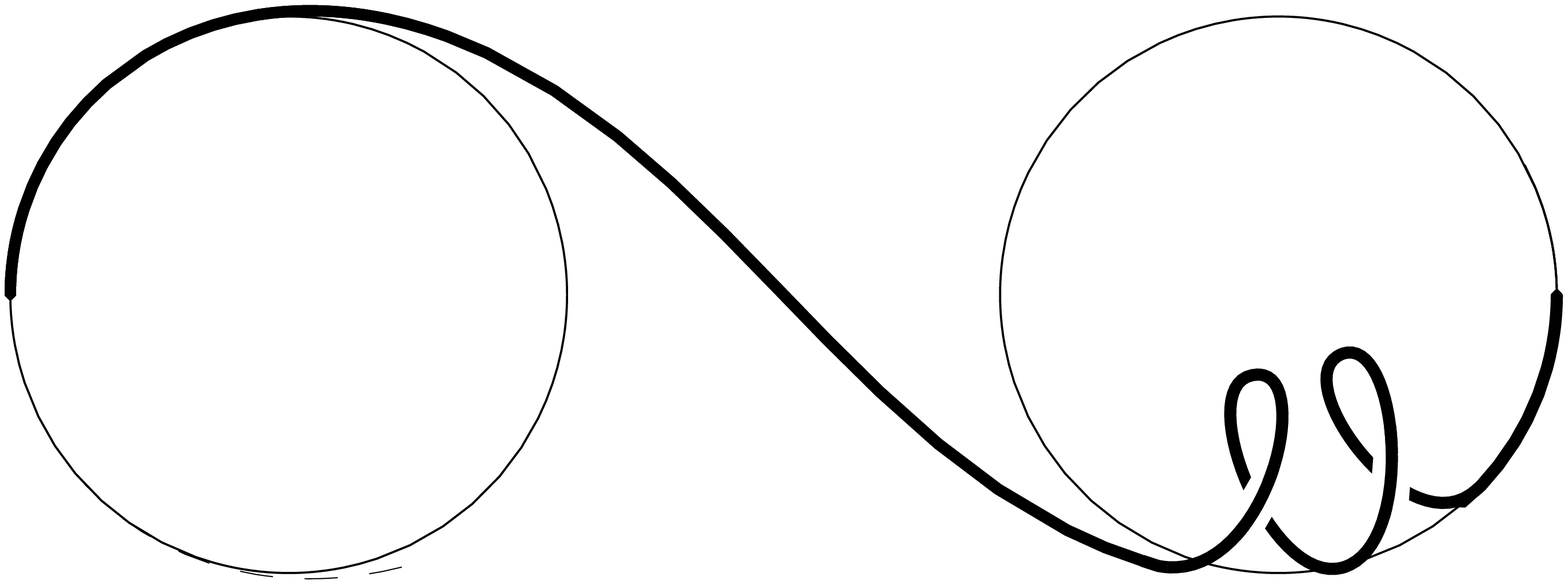}
\end{center}
\caption{$\sigma_{2}$ and $\sigma_{-2}$. }
\label{cate34}
\end{figure}

\medskip
\noindent
{\bf (The looped $S$-arc)}
Let $m$ be an integer, we attach $|m|$ loops 
to the $S$-arc $\sigma_0$, which
lies in the $xy$-plane as in Figure \ref{cate23} right.
Now, we slightly deform it as a space curve 
so that it has no self-intersection.
Figure \ref{cate34} left (resp. right)
indicates this new curve,
which is called the {\it $m$-looped $S$-arc}.
We denote it by $\sigma_m$.
Consequently, the $m$-looped $S$-arc is embedded,  
lies almost in the $xy$-plane, and 
has exactly one inflection point which
is just the original inflection point of the
lemniscate.

\begin{figure}[htb]
\begin{center}
 \includegraphics[height=4.0cm]{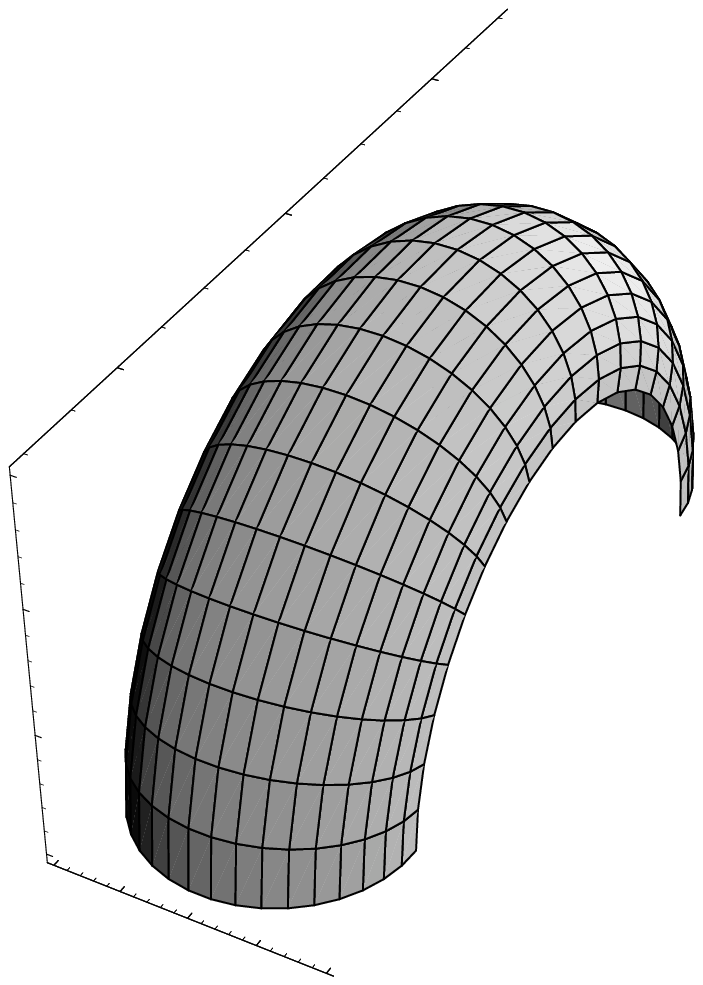}
\qquad
 \includegraphics[width=6.5cm]{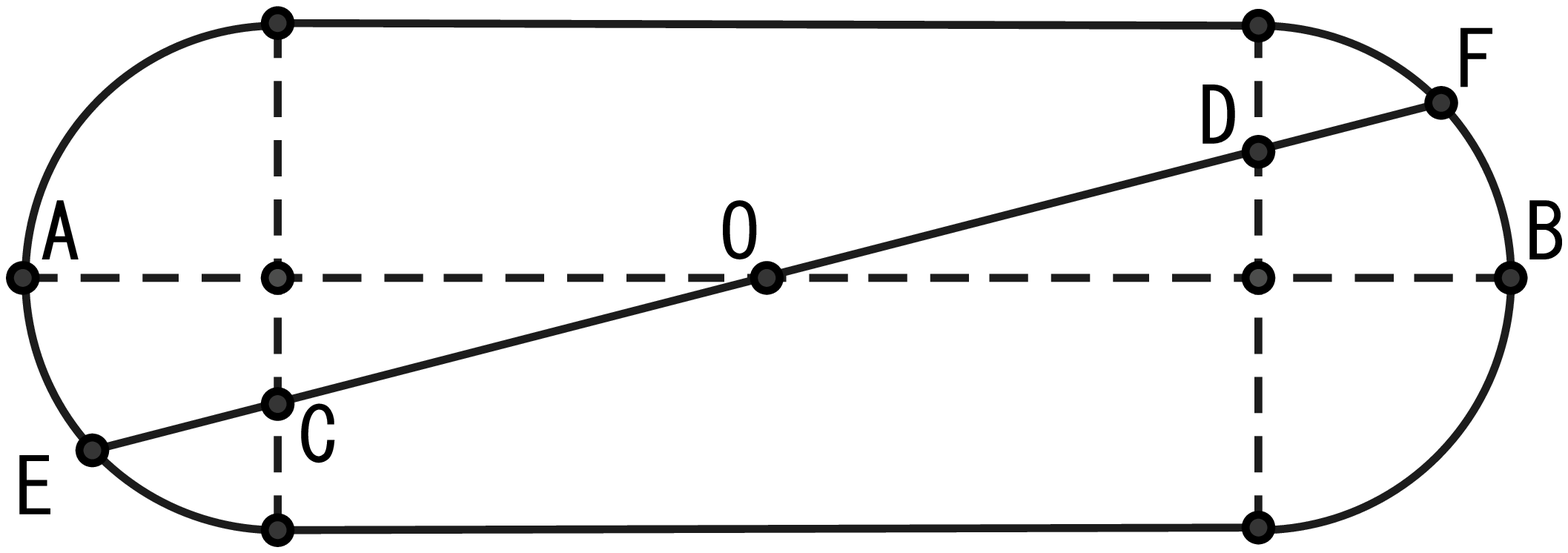}
\end{center}
\caption{The image of $f$ and $\Omega$.}
\label{torus}
\end{figure}

\medskip
\noindent
{\bf (The bridge arc on a torus)}
We set (cf. \eqref{eq:ab})
$$
a=\frac{2}{3},\qquad b=\frac{1}{3}
$$
and 
$$
f(u,v) := \pmt{(a + b \cos v)\cos u\\ b \sin v \\ (a + b \cos v)\sin u \\
  }
\qquad (|u|,|v|<\frac{\pi}2),
$$
which gives  an immersion into
the subset on a half-torus with positive 
Gaussian curvature as in Figure \ref{torus} left.
Then the two osculating circles at $t=\pi,2\pi$ (with radius $b$)
of the S-arc or the
looped S-arc (in $xy$-plane) lies in this torus.

Let 
$
\Pi:f([-\frac{\pi}2,\frac{\pi}2]\times [-\frac{\pi}2,\frac{\pi}2])
\to \R^2
$
be the projection into the $xy$-plane. Then the map $\Pi$ 
is injective, and the inverse map is given by
$$
\Pi^{-1}:\Omega\ni \pmt{x\\ y} \to 
\pmt{x\\ y\\ \left((a+\sqrt{b^2-y^2})^2-x^2\right)^{1/2}}\in \R^3,
$$
where $\Omega$ is the closed domain in the
$xy$-plane given by
$$
\Omega:=\{ |x|\le a, |y|\le b\}\cup
\{ (x-a)^2+y^2\le b^2\}\cup \{(x+a)^2+y^2\le b^2\}.
$$

\begin{figure}[htb]
\begin{center}
 \includegraphics[width=7.0cm]{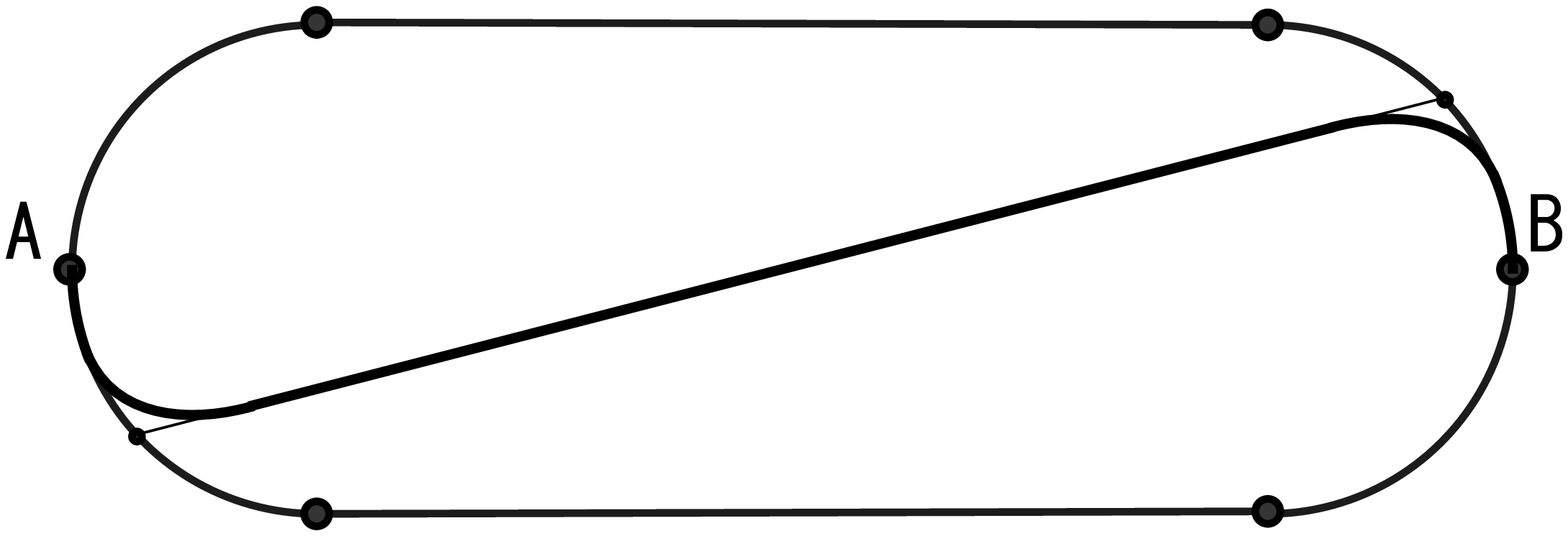}
\quad
 \includegraphics[width=7.0cm]{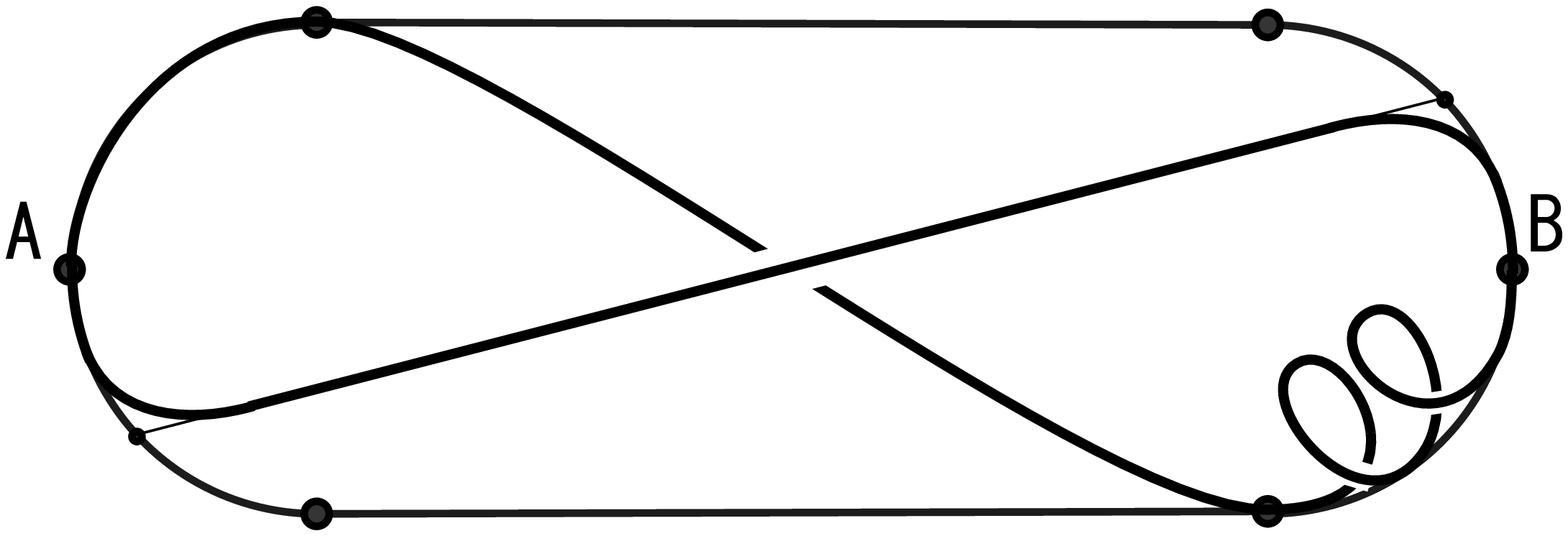}
\end{center}
\caption{$\tau_0$ and the top view of 
$c_0$.}
\label{torus-proj}
\end{figure}

We take the midpoints $\pt{A},\pt{B}$ on the circular parts on 
the boundary of $\Omega$.
Let $\pt{O}$ be the mid-point of $\pt{AB}$ which gives the 
center of gravity of $\Omega$.
Take two points $\pt{C},\pt{D}$ on $\Omega$ as in 
Figure \ref{torus} (right)
so that they bisect the radius of the circles of radius $b$.
Let $\pt{E},\pt{F}$ be the points 
where $\pt{CD}$ meets the 
boundary of $\Omega$.
We round the corner of the planar arc
$$
\overset{\frown}{\pt{AE}}\cup \pt{EF} \cup \overset{\frown}{\pt{FB}},
$$
and then we get a $C^\infty$-regular arc $\tau_0$ 
as in Figure \ref{torus-proj} (left).
The inverse image 
$$
\hat \tau_0:=\Pi^{-1}(\tau_0)
$$ 
on the torus is called the {\it bridge arc}.

\begin{lemma}\label{lem:tbridge}
Let $\hat \tau_0(t)$ $(0\le t\le 1)$ be the
the {\it bridge arc}. Then it has no inflections.
Moreover, it holds that
$$
\op{Tw}_{\hat \tau_0}(D^\perp)-
\op{Tw}_{\hat \tau_0}(\mb e_3^\perp)=\pi,
$$
where $D(t)$ is the Darboux vector field and 
$\mb e_3:=(0,0,1)$.
\end{lemma}
\begin{proof}
Let $\mb b(t)$ be the unit bi-normal vector of
$\hat \tau_0(t)$.
Let $\theta(t)$ be the smooth function 
which gives the leftward angle of 
$\mb b(t)$ from $\mb e_3^\perp$.
Like as the proof of Proposition \ref{prop:2-1},
we can see that $D^\perp_+=\mb b$.
Then we have
$$
\op{Tw}_{\hat \tau_0}(D^\perp)-
\op{Tw}_{\hat \tau_0}(\mb e_3^\perp)=\theta(1)-\theta(0).
$$
Let $\mb t(t)$ be the unit tangent vector
of $\tilde \sigma_+$ as a space curve.
Then by definition of $\tilde \sigma_+$,
we have 
$$
\mb t(0)=\mb t(1),\qquad \mb n(0)=-\mb n(1)
$$
which yield
\begin{equation}\label{eq:b12}
\mb b(0)=\mb t(0)\times \mb n(0)=
-\mb t(1)\times \mb n(1)=-\mb b(1).
\end{equation}
Since $\hat \tau_0(t)$ is planar near $t=0,1$,
$\mb b(t)$ is proportional to  $\mb e_3^\perp$.
Thus we have
$$
\theta(1)-\theta(0)\equiv \pi \,\, \mod 2\pi \Z.
$$
On the other hand,
the bridge arc $\hat \tau_0(t)$ ($0\le t\le 1$)
is symmetric with respect to the
plane containing the line $\pt{EF}$ which is perpendicular to
$xy$-plane.
Moreover,
the  bridge arc near the the mid point $\Pi^{-1}(O)$ 
is planar, and the 
$\mb b(t)$ is perpendicular to the plane.
Using these facts,
one can easily check that $\theta(t)\ge 0$, and
$$
\theta(1)-\theta(0)= \pi,
$$
which proves the assertion.
\end{proof}

Consider,
the union of the $m$-looped $S$-arc (a planar part)
and the bridge arc (a non-planar arc)
$$
(\mbox{Image of } \sigma_{m})
\cup \hat \tau_0,
$$
which gives a closed $C^\infty$-space curve.
We denote by
$c_0(t)=(x(t),y(t),z(t))$ ($|t|\le \pi$)
one of its parametrization.
Since $\hat \tau_0$ has no inflection points,
$c_0(t)$ is a closed embedded 
$C^\infty$-regular space curve with
a generic inflection point, which corresponds to
the inflection point of the original lemniscate.
Figure 11(right) shows the top view of $c_0$.

Without loss of generality,
we may assume that $t=0$ is the inflection point.
Let $D(t)$ be the normalized Darboux vector field along $c_0(t)$.
By Lemma \ref{lem:RR}, $F_{c_0,D}$ gives a 
rectifying unknotted $C^\infty$-M\"obius developable.
Moreover, by \eqref{eq:index} and Lemma \ref{lem:tbridge}
we can easily see that its M\"obius twisting number $2m-1$.
Since $m$ is arbitrary, its M\"obius twisting number can 
be adjusted arbitrarily.
Since the $S$-arc is planar, $c_0=(x(t),y(t),z(t))$ 
satisfies
$$
\dot z(0)=\ddot z(0)=\dddot z(0)=z^{(4)}(0)=0.
$$
On the other hand, rotating $F_{c_0,D}$ 
with respect to the $z$-axis,
we may assume
$$
\dot x(0)\ne 0,\quad \dot y(0)=\ddot y(0)=0,
\quad \dddot y(0)\ne 0,
$$
that is, $c_0(t)$ satisfies (1)-(3) of 
Proposition \ref{thm:admissible}.

\begin{figure}[htb]
\begin{center}
 \includegraphics[width=7.0cm]{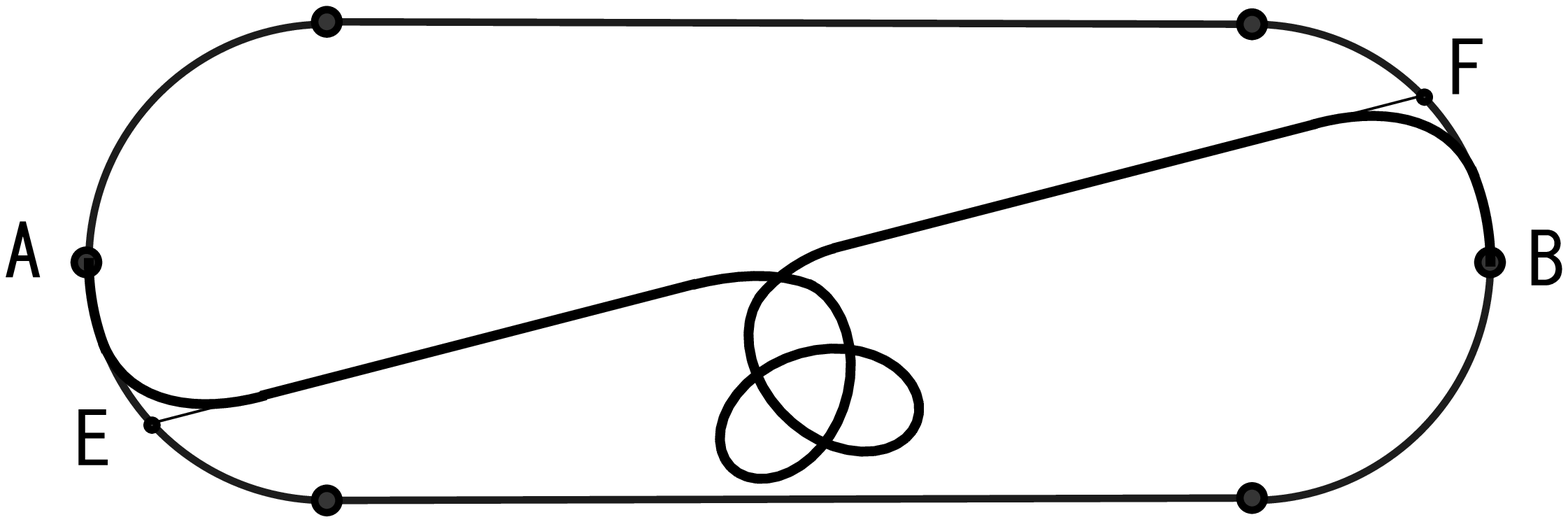}
\quad
 \includegraphics[width=7.0cm]{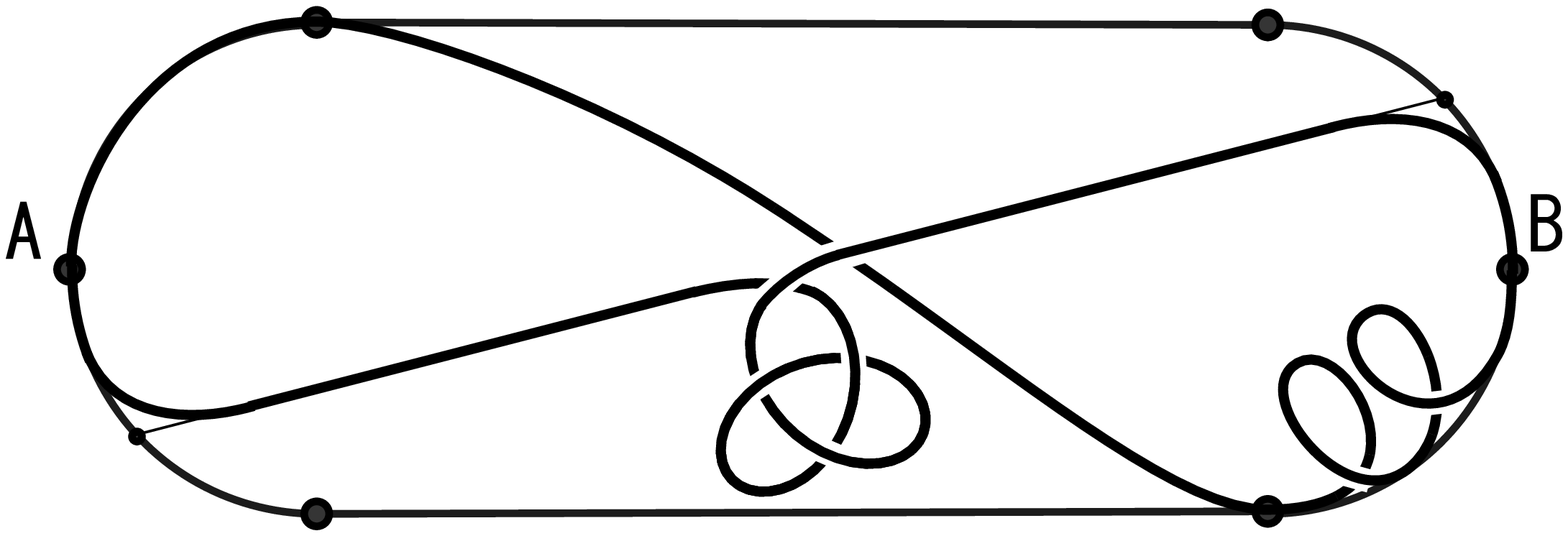}
\end{center}
\caption{$\tau_K$ and the top view of $c_K$}
\label{torus-proj2}
\end{figure}

Next, we construct a knotted
rectifying 
$C^\infty$ M\"obius developable.
Let $K$ be an arbitrary knot diagram.
Without loss of generality, we may assume that
$K$ lies in the subdomain $\{|x|<a,|y|<b\}$ of $\Omega$
and moreover that
$K$ lies in the lower half plane with respect to the line
$EF$. 
Like as in Figure 
\ref{torus-proj2} (left),
we connect $EF$ and the diagram $K$ , 
and denote it by $\tau_K$.
Since $\tau_K$ has self-intersections,
so is the inverse image $\Pi^{-1}(\tau_K)$.
By a small perturbation
near the each crossing (according to its up-down
status via the knot diagram $K$), 
we get an embedded
arc $\hat \tau_K$ 
on the torus, which  is called the {\it $K$-bridge arc}.
Since the Gaussian curvature on the torus
on $\Pi^{-1}(\tau_K)$ is positive,
the $K$-bridge arc $\hat \tau_K$ has no inflection points
as a space curve.
Let $c_K(t)=(x(t),y(t),z(t))$ ($|t|\le \pi$)
be the regular space curve which gives a
parametrization of 
the union of the $m$-looped $S$-arc and the
bridge arc 
$$
(\mbox{Image of } \sigma_{m})\cup \hat \tau_K.
$$
Then by the definition of $\hat \tau_K$,
$c_K$ is isotopic to the knot corresponding to $K$.
Moreover, $c_K$ gives an embedded 
$C^\infty$-regular space curve with
a generic inflection point.
Figure \ref{torus-proj2} right is the top view of $c_K$.
Without loss of generality we may assume that  $t=0$
is the inflection point.
Since the $S$-arc is planar, it satisfies
$$
\dot z(0)=\ddot z(0)=\dddot z(0)=z^{(4)}(0)=0,
$$
and we may assume
$$
\dot x(0)\ne 0,\quad 
\dot y(0)=\ddot y(0)=0,\quad \dddot y(0)\ne 0,
$$
like as in the unknotted case.
Consequently, $c_K$ satisfies (1)-(3) of 
Proposition \ref{thm:admissible}.
Taking the normalized Darboux vector field of 
$\hat \tau_K$, we get a rectifying 
M\"obius developable whose centerline is $\hat \tau_K$.
Since $m$ is arbitrary, its M\"obius twisting number can 
be adjusted arbitrarily.
Now we have just proved Proposition \ref{thm:admissible}.

\medskip
\noindent
Next we prove the following assertion:

\begin{proposition}\label{thm:final}
Let $\gamma(t)$ $(|t|\le \pi)$ be
a centerline of rectifying $C^\infty$-M\"obius developable $f$
satisfying the conditions $(1)$-$(3)$ in
Proposition \ref{thm:admissible}.
Then there exists a  family $\{\Gamma_n(t)\}$ 
$(|t|\le \pi)$ of
real analytic space curves such that
\begin{enumerate}
\item[{\rm (a)}] Each $\Gamma_n$
also satisfies conditions $(1)$-$(3)$ in
Proposition \ref{thm:admissible}.
\item[{\rm (b)}] $\{\Gamma_n\}_{n=1,2,\cdots}$ converges to $\gamma$
uniformly. Moreover, family of the $k$-th derivatives $(k=1,2,3,...)$
$\{\Gamma^{(k)}_n\}$ converges to $\gamma^{(k)}$
$C^\infty$-uniformly.
\end{enumerate}
In particular, the rectifying developable associated with
 $\Gamma_n$ converges $f$ uniformly.
\end{proposition}

\begin{proof}
We set
$$
\gamma(t)=(x(t),y(t),z(t))\qquad (| t |\le \pi).
$$
Consider a Fourier expansion of $\gamma(t)$
$$
\gamma(t)=a_0+\sum_{n=1}^\infty
\biggl( 
a_n \cos(nt) + b_n \sin(nt)
\biggr),
$$
and let
$$\gamma_n(t)\left (=(x_n(t),y_n(t),z_n(t))\right)
:=a_0+\sum_{j=1}^n
\biggl( 
a_j \cos(jt) + b_j \sin(jt)
\biggr)
\qquad (n=1,2,3,...)
$$
be the $n$ th approximation of $\gamma(t)$.
Then $\{\gamma_n\}$ is real analytic and
$C^\infty$-uniformly converges to $\gamma$.
Now we set
$$
X_n(t):=x_n(t),\qquad
Y_n(t):=y_n(t)-\dot y_n(0)\sin t+\ddot y_n(0)\cos t.
$$
Then they are real analytic
and satisfy
$$
\dot Y_n(0)=\ddot Y_n(0)=0.
$$
On the other hand, we have
$$
\dddot Y_n(t):=\dddot y_n(t)+\dot y_n(0)
\cos t+\ddot y_n(0)\sin t.
$$
Since
$$
\lim_{n\to 0}\dot y_n(0)
=\lim_{n\to 0} \ddot y_n(0)= 0,\quad
\lim_{n\to 0}\dddot y_n(0)=\dddot y(0)\ne 0,
$$
we have
$$
\dddot Y_n(0)\ne 0
$$
for sufficiently large $n$.
Next we set
\begin{align*}
Z_n(t):=& z_n(t)+\frac{4\ddot z_n(0)+z^{(4)}_n(0)}{3}\sin t-
\frac{4\dot z_n(0)+\dddot z_n(0)}{3}\cos t \\
&\qquad \quad \qquad -\frac{\ddot z_n(0)+z^{(4)}_n(0)}{12}\sin(2t)
+\frac{\dot z_n(0)+\dddot z_n(0)}{6}\cos(2t).
\end{align*}
Then it satisfies
$$
\dot Z_n(0)=\ddot Z_n(0)=\dddot Z_n(0)=Z^{(4)}_n(0)=0.
$$
If we set
$$
\Gamma_n(t)=(X_n(t),Y_n(t),Z_n(t)),
$$
then it satisfies (2) and (3) of 
Proposition \ref{thm:admissible}.
Moreover, we have
\begin{equation}\label{eq:condA}
\lim_{n\to 0}\dot y_n(0)=\lim_{n\to 0} \ddot y_n(0)= 
\lim_{n\to 0}\dot z_n(0)=\lim_{n\to 0} \ddot z_n(0)= 
\lim_{n\to 0}\dddot z_n(0)=\lim_{n\to 0} z^{(4)}_n(0)= 0.
\end{equation}
Since 
$\gamma_n$  converges $C^\infty$-uniformly to $\gamma$,
so does $\Gamma_n$ because of \eqref{eq:condA}.

Next we show that
$\Gamma_n(t)$ ($t\ne 0$) has no inflection point.
It can be checked by a straight-forward calculation
that
$
\dot \Gamma_n\times \ddot \Gamma_n
$
converges to $\dot \gamma\times \ddot \gamma$ on $[-\pi,\pi]$
uniformly.
Thus for any $\varepsilon>0$,
there exists a positive integer $N$
such that 
$\Gamma_n(t)$ ($n\ge N$) has no inflection point
for $|t|\ge \epsilon$.
So it is sufficient to prove that 
there exists $\epsilon>0$ such that
$\dot\Gamma_n(t)\times \ddot \Gamma_n(t)$
($|t|<\epsilon$) vanishes only at $t=0$:
The third component of the binormal vector
$$
(\beta_1,\beta_2,\beta_3):=\dot \Gamma_n(t)
\times \ddot \Gamma_n(t)
$$ 
is given by
$$
\beta_3(t)
=\ddot x_n(t)(-\dot y(t)+\dot y(0)\cos t+\ddot y(0)\sin t )
+
\dot x(t)(\ddot y(t)+\dot y(0)\sin t-\ddot y(0)\cos t ).
$$
Since
\begin{align*}
\sin t&=t +o(t),\qquad \cos t=1+o(t), \\
\dot x(t)&=\dot x(0)+o(t),\\
\dot y(t)&=\dot y(0)+t \ddot y(0)+o(t),\\
\ddot y(t)&=\ddot y(0)+t \dddot y(0)+o(t),
\end{align*}
we have
$$
\beta_3(t)=\dot x_n(0)(\dddot y_n(0)+\dot y_n(0))t+o(t^2).
$$
Here $o(t)$ and $o(t^2)$ are the higher order terms
than $t$ and $t^2$ at $t=0$, respectively.
Since
$$
\lim_{n\to \infty}\dot x_n(0)=\dot x(0)\ne 0,\quad
\lim_{n\to \infty}\dot y_n(0)=0,\quad
\lim_{n\to \infty}\dddot y_n(0)=\dddot y(0)\ne 0,
$$
we can conclude that
$\dot \Gamma_n(t)\times \ddot \Gamma_n(t)$
does not vanish for sufficiently small $t\ne 0$
and for sufficiently large $n$.

Finally, we show that
the rectifying developable associated with
 $\Gamma_n$ converges $f$ uniformly.
Then the Darboux vector field $D_n(t)$ of 
$\Gamma_n(t)$ has the following expression
$$
D_n(t)=\frac{\tau_n(t)}{\kappa_n(t)}\mb t_n(t) +
\mb b_n(t)
$$
for $t\ne 0$, where $\mb t_n$, 
$\mb b_n$, $\kappa_n$ and $\tau_n$
are unit tangent vector, the unit bi-normal 
vector, the curvature and the torsion respectively.

Since $\Gamma_n(t)$ is real analytic and $t=0$
is a generic inflection point, 
there exists a real analytic $\R^3$-valued function
$\mb c_n(t)$ such that $\mb c_n(0)\ne 0$ and
$$
\dot \Gamma_n(t)\times \ddot \Gamma_n(t)=t \mb c_n(t).  
$$
Then
$$
\mb b_n(t)=\frac{\mb c_n(t)}{|\mb c_n(t)|}
$$
gives a smooth parametrization of unit bi-normal vector
of $\Gamma_n(t)$ near $t=0$.
On the other hand, Let $M$ be the order of torsion at $t=0$.
Since $\Gamma_n(t)$ satisfies (1)-(3) of Proposition \ref{thm:admissible}, 
we have $M\ge 3$.
Since $\Gamma_n(t)$ is real analytic, 
there exists a real analytic $\R^3$-valued function
$T_n(t)$ such that
$$
\op{det}(\dot \Gamma_n(t),
\ddot\Gamma_n(t),
\dddot \Gamma_n(t))
=t^3 T_n(t).  
$$
Thus we have
$$
\frac{\tau_n(t)}{\kappa_n(t)}
=\frac{\op{det}(\dot \Gamma_n(t),\ddot\Gamma_n(t),
\dddot \Gamma_n(t))}{|\dot \Gamma_n(t)\times \ddot \Gamma_n(t)
|^3}=
\frac{T_n(t)}{|\mb c_n(t)|^3}.
$$
Since $\Gamma_n(t)$ converges to $\gamma(t)$ $C^\infty$-uniformly, 
The normalized Darboux vector field $D_n(t)$ also converges 
uniformly to that of $\gamma(t)$.
\end{proof}

\medskip
\noindent
({\it Proof of Theorem B.})
There exists an embedded rectifying $C^\infty$ M\"obius developable $F$
with an arbitrary isotopy type such that its centerline
$$
\gamma(t)=(x(t),y(t),z(t))\qquad (|t|\le \pi )
$$
as a $2\pi$-periodic embedded space curve satisfying
satisfying (1)-(3) of Proposition \ref{thm:admissible}.
By Proposition \ref{thm:final} and Corollary \ref{cor:generic},
there exists a sequence $\{F_n\}$ of
rectifying $C^\omega$ M\"obius developable uniformly
converges to $F$.
Then $F_n$ is the same isotopy type as $F$ if $n$ is sufficiently
large. \rightline (q.e.d.)

\begin{acknowledgements}
The authors thank Ryushi Goto 
for fruitful conversations on this subject.
The authors also thank 
Wayne Rossman for a careful reading of the first draft 
and for valuable comments.
\end{acknowledgements}


\end{document}